\documentclass[11pt,reqno]{amsart}
\usepackage{amsthm}
\usepackage{amssymb}
\usepackage{latexsym}
\usepackage{multicol}
\usepackage{verbatim}
\usepackage{amscd}
\usepackage{hyperref}
\usepackage{amsmath, amscd}
\usepackage[all]{xy}
\usepackage{color}
\advance\textwidth by 1.2in \advance\oddsidemargin by -.6in \advance\evensidemargin by -.6in

\newtheorem*{cor}{Corollary}
\newtheorem*{lem}{Lemma}
\newtheorem*{prop}{Proposition}

\theoremstyle{definition}
\newtheorem*{defn}{Definition}
\theoremstyle{definition}
\newtheorem*{thm}{Theorem}
\newtheorem{Theorem}{Theorem}
\newtheorem*{conj}{Conjecture}

\newenvironment{pf}{\proof}{\endproof}
\newcounter{cnt}
\newenvironment{enumerit}{\begin{list}{{\hfill\rm(\roman{cnt})\hfill}}{%
\settowidth{\labelwidth}{{\rm(iv)}}\leftmargin=\labelwidth%
\advance\leftmargin by \labelsep\rightmargin=0pt\usecounter{cnt}}}{\end{list}} \makeatletter
\def\mydggeometry{\makeatletter\dg@YGRID=1\dg@XGRID=20\unitlength=0.003pt\makeatother}
\makeatother \theoremstyle{remark}

\numberwithin{equation}{section}

\let\bwdg\bigwedge
\def\bigwedge{{\textstyle\bwdg}}

\begin{document}
\newcommand{\thmref}[1]{Theorem~\ref{#1}}
\newcommand{\secref}[1]{Section~\ref{#1}}
\newcommand{\lemref}[1]{Lemma~\ref{#1}}
\newcommand{\propref}[1]{Proposition~\ref{#1}}
\newcommand{\corref}[1]{Corollary~\ref{#1}}
\newcommand{\remref}[1]{Remark~\ref{#1}}
\newcommand{\defref}[1]{Definition~\ref{#1}}
\newcommand{\er}[1]{(\ref{#1})}
\newcommand{\id}{\operatorname{id}}
\newcommand{\ord}{\operatorname{\emph{ord}}}
\newcommand{\sgn}{\operatorname{sgn}}
\newcommand{\wt}{\operatorname{wt}}
\newcommand{\tensor}{\otimes}
\newcommand{\from}{\leftarrow}
\newcommand{\nc}{\newcommand}
\newcommand{\rnc}{\renewcommand}
\newcommand{\dist}{\operatorname{dist}}
\newcommand{\qbinom}[2]{\genfrac[]{0pt}0{#1}{#2}}
\nc{\cal}{\mathcal} \nc{\goth}{\mathfrak} \rnc{\bold}{\mathbf}
\renewcommand{\frak}{\mathfrak}
\newcommand{\supp}{\operatorname{supp}}
\newcommand{\Irr}{\operatorname{Irr}}
\renewcommand{\Bbb}{\mathbb}
\nc\bomega{{\mbox{\boldmath $\omega$}}} \nc\bpsi{{\mbox{\boldmath $\Psi$}}}
 \nc\balpha{{\mbox{\boldmath $\alpha$}}}
 \nc\bpi{{\mbox{\boldmath $\pi$}}}
\nc\bsigma{{\mbox{\boldmath $\sigma$}}} \nc\bcN{{\mbox{\boldmath $\cal{N}$}}} \nc\bcm{{\mbox{\boldmath $\cal{M}$}}} \nc\bLambda{{\mbox{\boldmath
$\Lambda$}}}

\newcommand{\lie}[1]{\mathfrak{#1}}
\makeatletter
\def\section{\def\@secnumfont{\mdseries}\@startsection{section}{1}%
  \z@{.7\linespacing\@plus\linespacing}{.5\linespacing}%
  {\normalfont\scshape\centering}}
\def\subsection{\def\@secnumfont{\bfseries}\@startsection{subsection}{2}%
  {\parindent}{.5\linespacing\@plus.7\linespacing}{-.5em}%
  {\normalfont\bfseries}}
\makeatother
\def\subl#1{\subsection{}\label{#1}}
 \nc{\Hom}{\operatorname{Hom}}
  \nc{\mode}{\operatorname{mod}}
\nc{\End}{\operatorname{End}} \nc{\wh}[1]{\widehat{#1}} \nc{\Ext}{\operatorname{Ext}} \nc{\ch}{\operatorname{ch}} \nc{\ev}{\operatorname{ev}}
\nc{\symm}{\operatorname{symm}}
\nc{\Ob}{\operatorname{Ob}} \nc{\soc}{\operatorname{soc}} \nc{\rad}{\operatorname{rad}} \nc{\head}{\operatorname{head}}
\def\Im{\operatorname{Im}}
\def\gr{\operatorname{gr}}
\def\mult{\operatorname{mult}}
\def\Max{\operatorname{Max}}
\def\ann{\operatorname{Ann}}
\def\annh{\operatorname{Ann_{\bu(\lie h[t])}}}
\def\sym{\operatorname{sym}}
\def\Res{\operatorname{\br^\lambda_A}}
\def\und{\underline}
\def\Lietg{$A_k(\lie{g})(\bsigma,r)$}

 \nc{\Cal}{\cal} \nc{\Xp}[1]{X^+(#1)} \nc{\Xm}[1]{X^-(#1)}
\nc{\on}{\operatorname} \nc{\Z}{{\bold Z}} \nc{\J}{{\cal J}} \nc{\C}{{\bold C}} \nc{\Q}{{\bold Q}}
\renewcommand{\P}{{\cal P}}
\nc{\N}{{\Bbb N}} \nc\boa{\bold a} \nc\bob{\bold b} \nc\boc{\bold c} \nc\bod{\bold d} \nc\boe{\bold e} \nc\bof{\bold f} \nc\bog{\bold g}
\nc\boh{\bold h} \nc\boi{\bold i} \nc\boj{\bold j} \nc\bok{\bold k} \nc\bol{\bold l} \nc\bom{\bold m} \nc\bon{\bold n} \nc\boo{\bold o}
\nc\bop{\bold p} \nc\boq{\bold q} \nc\bor{\bold r} \nc\bos{\bold s} \nc\boT{\bold t} \nc\boF{\bold F} \nc\bou{\bold u} \nc\bov{\bold v}
\nc\bow{\bold w} \nc\boz{\bold z} \nc\boy{\bold y} \nc\ba{\bold A} \nc\bb{\bold B} \nc\bc{\bold C} \nc\bd{\bold D} \nc\be{\bold E} \nc\bg{\bold
G} \nc\bh{\bold H} \nc\bi{\bold I} \nc\bj{\bold J} \nc\bk{\bold K} \nc\bl{\bold L} \nc\bm{\bold M} \nc\bn{\bold N} \nc\bo{\bold O} \nc\bp{\bold
P} \nc\bq{\bold Q} \nc\br{\bold R} \nc\bs{\bold S} \nc\bt{\bold T} \nc\bu{\bold U} \nc\bv{\bold V} \nc\bw{\bold W} \nc\bz{\bold Z} \nc\bx{\bold
x} \nc\KR{\bold{KR}} \nc\rk{\text{rk}} \nc\het{\text{ht }}
\nc\loc{\rm{loc }}

\nc\toa{\tilde a} \nc\tob{\tilde b} \nc\toc{\tilde c} \nc\tod{\tilde d} \nc\toe{\tilde e} \nc\tof{\tilde f} \nc\tog{\tilde g} \nc\toh{\tilde h}
\nc\toi{\tilde i} \nc\toj{\tilde j} \nc\tok{\tilde k} \nc\tol{\tilde l} \nc\tom{\tilde m} \nc\ton{\tilde n} \nc\too{\tilde o} \nc\toq{\tilde q}
\nc\tor{\tilde r} \nc\tos{\tilde s} \nc\toT{\tilde t} \nc\tou{\tilde u} \nc\tov{\tilde v} \nc\tow{\tilde w} \nc\toz{\tilde z}

\title{BGG reciprocity for current algebras}
\author{Matthew Bennett, Vyjayanthi Chari \and Nathan Manning}
\thanks{VC acknowledges support from the NSF, DMS-0901253}
\email{chari@math.ucr.edu }
\email{nmanning@math.ucr.edu}
\email{mbenn002@math.ucr.edu}

\begin{abstract} We study the category $\cal I_{\gr}$  of graded representations with finite--dimensional graded pieces for the current algebra $\lie g\otimes\bc[t]$ where $\lie g$ is a simple Lie algebra. This category has   many similarities with the category $\cal O$ of modules for $\lie g$ and  in this paper, we formulate and  study an  analogue of the famous BGG duality. We recall the definition of the projective and simple objects in $\cal I_{\gr}$ which are indexed by dominant integral weights. The role of the Verma modules is played by a family of modules called the global Weyl modules. We show that in the case when $\lie g$ is of type $\lie{sl}_2$,  the projective module admits a flag in which the successive quotients are finite direct sums of global Weyl modules. The multiplicity with which a particular Weyl module occurs in the flag  is determined by the multiplicity of a Jordan--Holder series for a closely associated family of modules, called the local Weyl modules. We conjecture that the result remains true for arbitrary simple Lie algebras. We also prove some combinatorial product--sum identities involving Kostka polynomials which arise as a consequence of our theorem.

\end{abstract}

\maketitle

\section*{Introduction} The category $\cal I_0$ of level zero  representations and its full subcategory $\cal I_{fin}$ of finite dimensional representations of affine and quantum affine algebras have been extensively studied in recent years.  The simple objects in these categories were classified in \cite{Cinv}, \cite{CPbanff}, \cite{CG1}. In the  case of affine Lie algebras, it was relatively easy to describe these simple objects \cite{CPnew}. In fact the classification of simple objects has been extended recently to the case of equivariant map algebras \cite{NSS}, \cite{NS}.

The irreducible objects in the  quantum situation, however, have proved to be very difficult to understand and an extensive list  of references can be found in \cite{CHer}. One method used to study the simple objects in the quantum case is to understand the $q=1$ limit of these representations, in which case they become indecomposable and generally reducible representations of the affine Lie algebra. In fact, in many situations, they remain indecomposable as representations of the maximal parabolic subalgebra of the affine algebra: namely of the Lie  algebra $\lie g[t]$ of polynomial maps from $\bc\to \lie g$, where $\lie g$ is the underlying simple Lie algebra.

The current  paper was motivated by an attempt to analyze the categories $\cal I_0$ and $\cal I_{fin}$ by understanding its homological properties and to see if one could apply the theory of highest weight categories of \cite{CPS}. We were also motivated by the classical Bernstein--Gelfand--Gelfand (BGG) result of \cite{BGG}; a more recent exposition can be found in \cite{Hu2}. To do this one needs an interval--finite poset which indexes three families of suitably related objects in the category:  the simple objects, projective covers of the simple objects and the standard modules.
The problem in our case, however, is that the categories $\cal I_0$ and $\cal I_{fin}$ do not have enough projectives and  the irreducible objects are indexed by continuous parameters. However, in \cite{CPweyl} a family of discretely indexed modules called the global Weyl modules were defined and they appeared to be a natural candidate for the standard modules. More evidence that the global Weyl modules were the correct objects was provided in the paper \cite{BCGM}, where the homomorphisms between global Weyl modules were studied. It became natural then, to look for suitable subcategories of $\cal I_0$ in which one could develop the theory of highest weight categories.
Thus, in \cite{CG1} it was seen  that  the category $\mathcal I_{\gr}$ whose objects are $\bz$--graded $\lie g[t]$ representations, with finite--dimensional graded components and graded morphisms had many nice properties. The simple objects in the category are  indexed by a discrete set and have projective covers. Moreover, the global Weyl modules lie in $\cal I_{\gr}$, as does a distinguished family of quotients of these modules, the so-called graded local Weyl modules. There are a number of reasons for preferring to work with the current algebra rather than the loop algebra; one of them, certainly, is the fact that there are no graded local Weyl modules in the case of the loop algebra. However, some of the results of this paper should also be true for the loop algebra.

We describe our results in some detail. If $P^+$ denotes the set of dominant integral weights of $\lie g$ and $\bz$ denotes the ring of integers, then it is easily seen \cite{CG1} that
the isomorphism classes of simple objects of $\mathcal I_{\gr}$ are indexed by $P^+\times\bz$.  For $(\lambda,r)\in P^+\times\bz$, the simple object $V(\lambda,r)$ of $\mathcal I_{\gr}$ is given in a natural way by fixing the irreducible finite--dimensional $\lie g$-module $V(\lambda)$ to lie in grade $r$ and assigning the appropriate $\lie g\otimes\bc[t]$-module structure. The projective cover $P(\lambda,r)$ of $V(\lambda,r)$ is  defined canonically by a standard construction of relative homological algebra, and the global Weyl module $W(\lambda,r)$ is  then seen to be a truncation of  $P(\lambda,r)$ obtained by requiring the set of $\lie g$--weights of $W(\lambda,r)$ to lie in the cone below $\lambda$.
We define the notion of a Weyl flag for an object $M$  of $\mathcal I_{\gr}$ in an obvious way: namely, it is  a descending filtration of  submodules of $M$, in which the successive quotients are  isomorphic to  global Weyl modules.  We show that the multiplicity with which a given Weyl module appears, is independent of the Weyl flag.

 Our main result is that when $\lie g$ is of type $\lie{sl}_2$, the object $P(\lambda,r)$ admits a Weyl flag. To formulate a BGG-type duality, one needs to show that the graded multiplicity of $W(\mu,s)$ is the same as the Jordan--Holder multiplicities $V(\lambda,r)$ in $W(\mu,s)$. But here one runs into a problem since the global Weyl modules are of infinite length. Our key observation is that this problem can be avoided if one works instead with the corresponding graded local Weyl module, which is finite--dimensional and hence has a Jordan--Holder series. This allows us to state and prove the correct analog of the BGG reciprocity between the Weyl filtration of the projective modules  and the Jordan--Holder series of the  graded local Weyl module.

We now   explain  the organization and structure of the paper.  Section \ref{prelim} introduces the notation and basic results, together with some general constructions. In Section \ref{mainresult}, we define the principal objects of study for arbitrary simple Lie algebras  and state the main result  of the paper for $\lie{sl}_2$. We also give a natural conjecture for  the general case.
Section \ref{characters} is devoted to recalling necessary results on local and global Weyl modules from  \cite{CL1},  \cite{CPweyl}, \cite{FoL}, \cite{Naoi}. We also recall the formulae of the graded characters of the local  Weyl modules from \cite{CL1}, \cite{Naoi} and deduce the graded character of the global Weyl modules.  As an application of our main theorem,   we then deduce an interesting combinatorial identity (see Section \ref{kostkasl2}) given in terms of the Kostka polynomials and compute the graded $\lie g$--character of the universal enveloping algebra of $\lie g\otimes t\bc[t]$. Section \ref{examples} simultaneously provides non--trivial examples of modules admitting a Weyl flag for an arbitrary simple Lie algebra  and also offers some evidence for our conjecture. The last two sections of the paper are devoted to proving the main theorem.

We conclude the introduction with the following remarks. Once BGG--duality is established in $\cal I_{\gr}$, one can then naturally find full subcategories which are highest weight. These ideas will be explored further in a future publication. It would also be interesting to see to what extent the results of this paper could be generalized to equivariant map algebras. The corresponding Weyl modules have been studied in the untwisted  cases in \cite{CFK} and also in \cite{FKKS} in the twisted case.

{\em Acknowledgements.} We thank Jacob Greenstein,  Sergey Loktev and David Rush  for helpful discussions.  It is a pleasure for the second and third authors to thank the organizers for their invitation to participate in the  trimester  \lq\lq On the interactions of Representation theory with Geometry and Combinatorics,\rq\rq at the Hausdorff Institute, Bonn, 2011.

\section{Preliminaries}\label{prelim} We establish the notation to be used in the rest of the paper and recall some standard results.
\subsection{} Throughout this paper we let   $\bc$  be the field of complex numbers and $\bz$ (resp. $\bz_+$, $\bn$) the set of integers (resp.  nonnegative, positive integers). Given a $\bz$--graded complex vector space $V=\oplus_{s\in\bz}V[s]$ with $\dim V[s]<\infty$ for all $s\in\bz$, let   $$\mathbb H(V)=\sum_{s\in\bz}\dim V[s]u^s\in\bc[[u,u^{-1}]],$$  be the Hilbert series in an indeterminate   $u$.

For a Lie algebra $\lie a$, we denote by $\bu(\lie a)$ the universal enveloping algebra of $\lie a$. If $\lie b$ and $\lie c$ are subalgebras of $\lie a$, with $\lie a=\lie b\oplus\lie c$ as vector spaces, then we have an isomorphism of vector spaces $$\bu(\lie a)\cong\bu(\lie b)\otimes\bu(\lie c).$$ The assignment $a\mapsto a\otimes 1+1\otimes a$, for $a\in\lie a$, defines a Hopf algebra structure on $\bu(\lie a)$ and we let $\bu(\lie a)_+$ be the augmentation ideal of $\bu(\lie a)$. We say that $\lie a$ is a $\bz$--graded Lie algebra if $$\lie a=\bigoplus_{s\in\bz}\lie a[s],\ \qquad [\lie a[r],\lie a[s]]\subset\lie a[r+s],\ \ r,s\in\bz,$$ and in this case the algebra $\bu(\lie a)$ inherits a $\bz$--grading.

{\em For the rest of the paper  we will drop the dependence on $\bz$, in other words when we refer to a graded vector space or graded Lie algebra and so on, we will mean a $\bz$--graded vector space, Lie algebra etc.}

A graded module for a  graded Lie algebra $\lie a$ is a graded complex vector space $V =\bigoplus_{s\in\bz}V[s]$ which admits a left action of $\lie a$ which is compatible with the grading, i.e, $(\lie a[r])V[s]\subset V[s+r]$ for all $ s,r\in\bz.$ Given two graded modules $V$, $W$ for $\lie a$ we say that a map $\pi: V\to W$ is a map of graded $\lie a$--modules if it is a map of $\lie a$--modules and  $\pi V[r]\subset W[r]$ for all $r\in\bz$.

We shall be interested in graded Lie algebras which are obtained as follows: given a  complex Lie algebra $\lie b$ and  an  indeterminate $t$, let  $\lie b[t]=\lie b\otimes \bc[t]$ be  the Lie algebra with commutator given by, $$[a\otimes f, b\otimes g]=[a,b]\otimes fg,\ \ a, b\in\lie b, \ \ f,g,\in\bc[t].$$  We identify without further comment the Lie algebra $\lie b$ with the subalgebra $\lie b\otimes 1$ of $\lie b[t]$. Clearly,  $\lie b[t]$ has a natural $\bz$--grading given by the powers of $t$ and this also induces a $\bz$--grading on $\bu(\lie b[t])$. Moreover, if $\dim \lie b<\infty$, then $\dim\bu(\lie b\otimes t\bc[t])[r]<\infty$ for all $r\in\bz$. Notice that  $$\bu(\lie b[t])[s]=0,\ \ \  s<0,\ \ \qquad  \bu(\lie b[t])[0]=\bu(\lie b).$$

\subsection{}
From now on, $\lie g$ will denote a finite--dimensional, simple Lie algebra of rank $n$ over the complex numbers, $\lie h$   a fixed Cartan subalgebra of $\lie g$ and $R$  the set of roots of $\lie g$ with respect to $\lie h$. Set   $I=\{1,\cdots ,n\}$  and choose $\{\alpha_i: i\in I\}$ a set of simple roots for $R$ and $\{\omega_i: i\in I\}$  a set  of fundamental weights.  Let $Q$ (resp. $Q^+$) be  the integer span (resp. the nonnegative integer span) of $\{\alpha_i: i\in I\}$ and similarly define  $P$ (resp. $P^+$) be the $\bz$ (resp. $\bz_+$) span of  $\{\omega_i: i\in I\}$. Define a  partial order  on $\lie h^*$ so that, given $\lambda,\mu\in\lie h^*$, we say that  $\lambda\ge\mu$ iff $ \lambda-\mu\in Q^+$.

 Set $R^+=R\cap Q^+$ and for $\alpha\in R^+$ let $\lie g_\alpha$ be the corresponding root space of $\lie g$ and  set $\lie n^\pm=\bigoplus_{\alpha\in R^+}\lie g_{\pm\alpha}$. Fix non--zero elements  $x^\pm_\alpha\in\lie g_{\pm\alpha}$, $h_\alpha\in\lie h$  such that $$[h_\alpha, x^\pm_\alpha]=\pm 2x^\pm_\alpha,\ \  [x^+_\alpha, x^-_\alpha]=h_\alpha,$$ and  for $i\in I$, we also write $x^\pm_i$, $h_i$ for $x^\pm_{\alpha_i}$, $h_{\alpha_i}$.
We have, $$\lie g\ =\ \lie n^-\oplus\lie h\oplus\lie n^+,\ \ \qquad \bu(\lie g)=\bu(\lie n^-)\otimes\bu(\lie h)\otimes\bu(\lie n^+). $$
  Moreover, the adjoint action of $\lie h$ on $\lie g$ induces a decomposition,$$ \bu(\lie g)=\bigoplus_{\eta\in Q}\bu(\lie g)_{\eta},
\qquad \bu(\lie n^\pm)=\bigoplus_{\eta\in Q^+}\bu(\lie n^\pm)_{\pm\eta},$$ where $
 \bu(\lie g)_\eta=\{g\in\bu(\lie g): [h,g]=\eta(h)g,\ \ h\in \lie h\},$ and $\bu(\lie n^\pm)_\eta=\bu(\lie g)_{\eta}\cap\bu(\lie n^\pm)$.

 We  also have $$\lie g[t]=\lie n^-[t]\oplus\lie h[t]\oplus\lie n^+[t],\ \ \qquad \bu(\lie g[t]) =\bu(\lie n^-[t])\otimes \bu(\lie h[t])\otimes\bu(\lie n^+[t]).$$ The adjoint action of $\lie h$ on $\lie g[t]$ is diagonalizable and for $\eta\in Q$ the subspaces $\bu(\lie g[t])_\eta$ etc.  are defined in the obvious way. Moreover these subspaces are graded and we have \[\bu(\lie g[t])_\eta[r]=\bu(\lie g[t])_\eta\cap\bu(\lie g[t])[r].\]

\subsection{}\label{locfin}
For any $\lie g$-module $M$ and $\mu\in\lie h^*$, set
\[M_\mu=\{m\in M\ :\ h.m=\mu(h)m,\quad h\in\lie h\}.\] We say that $M$ is a weight module
for $\lie g$ if \[M=\bigoplus_{\mu\in\lie h^\ast} M_\mu,\] and  set $\wt(M)=\left\{\mu\in\lie h^\ast\ :\ M_\mu\ne 0\right\}$ and for $v\in M_\mu$ we set $\wt v=\mu$. Any finite--dimensional $\lie g$--module is  a weight module.
   It is well-known that the set of isomorphism classes of irreducible finite--dimensional $\lie g$-modules is in bijective correspondence with $P^+$. For  $\lambda\in P^+$ we denote by $V(\lambda)$ the representative of the corresponding isomorphism class, which is generated by  a vector $v_\lambda$ with defining relations
  \[ \lie n^+.v_\lambda=0,\qquad h.v_\lambda=\lambda(h)v_\lambda,\qquad (x_i^-)^{\lambda(h_i)+1}.v_{\lambda}=0,\ \ \  h\in\lie h,\ \ i\in I.\]
 Moreover $\wt V(\lambda)\subset\lambda-Q^+\subset P$. The module $V(0)$ is the trivial module for $\lie g$ and we shall write it as $\bc$.

If $M$ is a finite--dimensional $\lie g$--module, we denote by $[M :V(\mu)]$ the multiplicity of $V(\mu)$ in a Jordan--Holder series for $M$. The character of $M$ is the  element of the integral group ring $\bz[P]$ given by,
\[ {\rm ch}_{\lie g} M=\sum_{\mu\in P}\dim_\bc M_\mu e(\mu),\] where $e(\mu)\in\bz[P]$ is the generator of the group ring corresponding to $\mu$. The following is standard.
\begin{lem}\label{ind} The characters   $\{\ch_{\lie g} V(\mu): \mu\in P^+\}$ are a  linearly independent subset of $\bz[P]$.\hfill\qedsymbol\end{lem}

We say that $M$ is locally finite--dimensional if it is a direct sum of finite--dimensional modules for $\lie g$, in which case $M$ is necessarily a weight module.
Using Weyl's theorem one knows that  a locally finite--dimensional $\lie g$-module $M$ is isomorphic to a direct sum of $V(\lambda)$, $\lambda\in P^+$ and  hence in particular, $\wt M\subset P$. Set \begin{equation}\label{ndecomp} M^{\lie n^+}=\{m\in M: \lie n^+ m=0,\}\ \qquad M^{\lie n^+}_\lambda=M^{\lie n^+}\cap M_\lambda\cong\Hom_{\lie g}(V(\lambda), M). \end{equation} We have $$ M_\lambda=(\lie n^-M\cap M_\lambda)\oplus M_\lambda^{\lie n^+},\qquad
M=\bigoplus_{\lambda\in P^+}\bu(\lie g)M_\lambda^{\lie n^+},$$ or equivalently, the isotypical component of $M$ corresponding to $\lambda\in P^+$ is the $\lie g$--submodule generated by $M_\lambda^{\lie n^+}$.

\subsection{}

Let $\mathcal I_{\rm gr}$ be the  category whose objects are graded $\lie g[t]$-modules $V$ with finite dimensional graded pieces and where the morphisms are  (degree zero)  maps  of graded $\lie g[t]$-modules. Clearly $$\lie g V[r]\subset V[r],\ \ \ r\in\bz.$$

For any $r\in\bz$ we let $\tau_r$ be the grading shift operator. Clearly $\cal I_{\gr}$ is closed under taking submodules, quotients and finite direct sums. Given $V,W\in\Ob\cal I_{\gr}$, the tensor product $V\otimes W$ is naturally graded by setting
$$(V\otimes W)[k]=\bigoplus_{r\in\bz} V[r]\otimes W[k-r].$$
However the graded components are no longer finite--dimensional and hence $\cal I_{\gr}$ is not closed under taking tensor products. However, we do have
\begin{lem} Suppose that $V,W\in\Ob\cal I_{\gr}$ are such that $V[r] = W[r]=0$ for all but finitely  many $r<0$ (resp. $r>0$). Then $V\otimes W\in\Ob\cal I_{\gr}$.\hfill\qedsymbol \end{lem}
In the course of this paper we shall only be considering tensor products of objects which satisfy the conditions of the lemma and we shall just use the lemma without comment.

\subsection{}
The graded character  of $V\in\Ob\mathcal I_{\rm gr}$ is the  element of the ring  $\bz[P][[u,u^{-1}]]$  given by
\[{\rm ch}_{\rm{gr}}V:=\sum_{r\in\bz} {\rm ch}_\lie{g}(V[r])u^{r}.\]   The following is  straightforward.
\begin{lem}\label{graiso} Let $V\in\Ob\cal I_{\gr}$ and $\mu\in P^+$.  Then $$V_\mu^{\lie n^+}= \bigoplus_{r\in\bz} V_\mu^{\lie n^+}[r],\ \qquad V_\mu^{\lie n^+}[r]= V_\mu^{\lie n^+}\cap V[r].$$ In particular,
\begin{equation}\ch_{\gr} V=\sum_{\mu\in P^+}\ch_{\lie g} V(\mu)\mathbb H(V_\mu^{\lie n^+}).\end{equation}
\hfill\qedsymbol
\end{lem}

\subsection{}\label{gammatrunc} Given any subset $\Gamma $ of $P^+$, let $\cal I^{\Gamma}_{\gr}$ be the full subcategory of $\cal I_{\gr}$ consisting of objects $V$ such that $\wt V\subset\Gamma-Q^+$ and if  $\Gamma=\{\lambda\}$, then  set  $\cal I^\lambda_{\rm gr}=\cal I^\Gamma_{\rm gr}$. Again $\cal I_{\gr}^\Gamma $ is closed under taking submodules, quotients and finite direct sums. 
For $V\in\Ob\cal I_{\gr}$, set \begin{equation}\label{defgv}\Gamma V= \sum_{\mu\notin\Gamma-Q^+}\!\!\bu(\lie g[t])V_\mu,\qquad V^\Gamma=V/\Gamma V.\end{equation}
Then $V^\Gamma$ is the maximal quotient of $V$ that lies in $\cal I_{\gr}^\Gamma$.
If $\pi:V\to U$ is a morphism of objects in $\cal I_{\gr}$, then $\pi(\Gamma V)\subset\pi(\Gamma U)$ and  we have an induced morphism $\pi^\Gamma: V^\Gamma\to U^\Gamma$ and hence we have a  functor  from  $\cal I_{\gr}\to\cal I_{\gr}^\Gamma$. In the case that $\Gamma=\{\lambda\}$, we shall write $V^\Gamma=V^\lambda$.

\begin{prop}\label{trunc} For all  $\Gamma'\subset \Gamma$  and  $V\in\Ob\cal I_{\gr}$ we have an isomorphism of objects in $I_{\gr}$, $$V^{\Gamma'}\cong (V^\Gamma)^{\Gamma'},\ \qquad \ \Gamma'V^\Gamma\cong\Gamma'V/\Gamma V.$$
\end{prop}
\begin{pf} Let $\varphi_\Gamma: V\to V^\Gamma$ and  $\varphi_{\Gamma'}: V\to V^{\Gamma'}$ be the canonical maps  of $\lie g[t]$--modules. Since $\Gamma V\subset\Gamma' V$ we have a natural induced map
$\varphi_{\Gamma,\Gamma'}: V^\Gamma\to V^{\Gamma'}$ such that
 $\varphi_{\Gamma,\Gamma'}.\varphi_\Gamma=\varphi_{\Gamma'}$. We first prove that $\ker\varphi_{\Gamma,\Gamma'}=\Gamma'V^\Gamma$ which establishes the first isomorphism of the proposition.
Thus if $v\in V^\Gamma$ is such that $\varphi_{\Gamma,\Gamma'}=0$, there exists $u\in V$ such that
$$\varphi_\Gamma(u)=v,\ \ \varphi_{\Gamma'}(u)=0,$$ in other words $u\in\Gamma'V$. Since $\varphi_\Gamma(\Gamma'V)\subset\Gamma'V^\Gamma$   we see that
$\ker\varphi_{\Gamma,\Gamma'}\subset\Gamma'V^\Gamma $. For the reverse inclusion notice that $$\wt V^{\Gamma'}\subset\Gamma'-Q^+\implies \varphi_{\Gamma,\Gamma'}(\Gamma'V^\Gamma)=0.$$

We now prove the second isomorphism. As we have seen above, the restriction of $\varphi_\Gamma$ to $\Gamma'V$ induces a homomorphism $\varphi:\Gamma'V\to \Gamma'V^\Gamma$ and since
$\ker\varphi=\ker\varphi_\Gamma$ we get an injective map $\Gamma'V/\Gamma V\to \Gamma'V^\Gamma$. We prove now that $\varphi$ is surjective which clearly completes the proof. Thus, let $v\in\Gamma'V^\Gamma$
and write   $$v=\sum_{s=1}^\ell g_sv_s,\ \ v_s\in (V^\Gamma)_{\mu_s},\ \ \mu_s\notin \Gamma'-Q^+.$$    Choose elements $u_s\in V_{\mu_s}$ such that
 $\varphi_\Gamma(u_s)=v_s$ and set $u=\sum_{s=1}^\ell g_su_s.$ Clearly $u\in\Gamma'V$ and $\varphi_\Gamma(u)=v$ as required.
\end{pf}

\section{The Main Result  }\label{mainresult}

In this section we recall the definitions  of the main objects of study: the simple modules, the indecomposable projective covers of simple modules and the Weyl modules. Our main result is the relation between these three families of modules.

\subsection{} Let $\ev_0: \lie g[t]\to\lie g$ be the homomorphism of Lie algebras $ev_0(x\otimes f)= f(0)x$. The kernel of this map is a graded ideal in $\lie g[t]$ and hence the pull--back $\ev_0^* V$ of  a $\lie g$--module $V$ is an object of $\cal I_{\gr}$. Set $\tau_r\ev_0^* V(\lambda)=V(\lambda,r)$ and we let $v_{\lambda,r}\in V(\lambda,r)$ be the element corresponding to $v_\lambda$.

The following is straightforward,  \cite[Proposition 1.3]{CG1}
and classifies the irreducible modules in $\mathcal I_{\rm gr}$.

\begin{lem}\label{simpleclass}
 A simple  object  of  $\mathcal I_{\rm gr}$  is isomorphic to $V(\mu,r)$ for a unique choice
of  $(\mu,r)\in P^+\times\bz$.
 The simple objects of $\mathcal I_{\rm gr}^\lambda$ are the $V(\mu,r)$, with $(\mu,r)\in P^+\times\bz$, $\mu\le \lambda$. The graded character of $V(\lambda,r)$ is $\ch_{\lie g}V(\lambda)u^r$.\hfill\qedsymbol\end{lem}

\subsection{}\label{projectives} For $\mu\in P^+$ set
\[ P(\mu,0)=\bu(\lie g[t])\otimes_{\bu(\lie g)} V(\mu),\qquad P(\mu,r)=\tau_rP(\mu,0),\] and let $p_{\mu,r}\in P(\mu,r)$ be the image under $\tau_r$ of the element $1\otimes v_\mu$. Note that $$P(\mu,r)[r]\cong_{\lie g} V(\mu),\ \qquad \ P(\mu,r)[s]=0,\ \ \ s<r.$$

The following proposition was proved in \cite[Proposition 2.1]{CG1} for $P(\mu,r)$.
\begin{prop}\label{pdefrel} Given $(\lambda,r)\in P^+\times\bz$  we have that $P(\lambda,r)$ is an indecomposable  projective object in $\mathcal I_{\rm gr}$.  It    is generated as a $\lie g[t]$-module by $p_{\lambda,r}$ with defining relations 
\[ \lie n^+.p_{\lambda,r}=0,\quad h.p_{\lambda,r}=\lambda(h)p_{\lambda,r},\quad (x_i^-)^{\lambda(h_i)+1}.p_{\lambda,r}=0,\quad\text{for}\   h\in\lie h,\  i\in I. \]
Moreover, considering $P(\lambda,r)[r]$ as a $\lie g$-- module, we have
$$[P(\lambda,r)[r]: V(\lambda,r)]=1,$$
and hence  $V(\lambda,r)$ is the unique irreducible quotient of $P(\lambda,r)$ and  $P(\lambda,r)$ is its projective cover in $\mathcal I_{\rm gr}$. \hfill\qedsymbol
\end{prop}
The following is immediate.
\begin{cor} For $(\lambda,r)\in P^+\times\bz$ and $\Gamma\subset P^+$, the object $P(\lambda,r)^\Gamma$ is projective in $\cal I_{\gr}^\Gamma$.\hfill\qedsymbol
\end{cor}
The graded character of the projective module is given in terms of tensor products of the symmetric powers of the adjoint representation of $\lie g$, \cite[Proposition 2.1]{CG1}. One of the applications of the results proved or conjectured in the current  paper is a much more explicit description of this character and this is discussed in the next section.

\subsection{}
For $\lambda\in P^+$ and $r\in\bz$, the graded global Weyl module $W(\lambda,r)$  is defined by
\[W(\lambda,r)=P(\lambda,r)^\lambda\] and we denote the image of $p_{\lambda,r}$ by $w_{\lambda,r}$.
 The following result proved in  \cite[Proposition 3.3]{CFK}  makes the connection with  the original definition of the global Weyl modules given in \cite{CPweyl}.
\begin{lem}\label{weyldefrel}
The global Weyl module $W(\lambda,r)$ is indecomposable, and is generated as a $\lie g[t]$-module by $w_{\lambda,r}$ with defining relations \[  \lie n^+[t].w_{\lambda,r}=0,\quad h.w_{\lambda,r}=\lambda(h)w_{\lambda,r},\quad (x_{i}^-)^{\lambda(h_{i})+1}.w_{\lambda,r}=0,\quad\text{for}\   h\in\lie h,\  i\in I. \] Moreover, $$\wt W(\lambda,r)\subset\lambda-Q^+,\ \ [W(\lambda,r)[r]: V(\lambda)]=1,$$ and $V(\lambda,r)$ is the unique irreducible quotient of $W(\lambda,r)$.

\hfill\qedsymbol\end{lem}

\subsection{}  The graded characters of the global Weyl modules are known and we shall review this in the next section. The following trivial consequence of Lemma \ref{weyldefrel} is used repeatedly. For $\mu\in P^+$ and $r\in\bz$, we have, \begin{equation}\label{triv} \ch_{\gr} W(\lambda,r)-\ch_{\lie g} V(\lambda)u^r \in u^{r+1}\bz[P][[u, u^{-1}]].\end{equation}
\begin{lem}\label{linind} Considered as elements of $\bc[P][[u, u^{-1}]]$, the  characters of the global Weyl modules are linearly independent over $\bc$.
\end{lem}

\begin{pf}
Suppose that we have  $$ \sum_{\ell=1}^{m} a_\ell {\rm{ch}}_{\rm{gr}} W(\lambda_\ell,r_\ell)=0, \ \ \qquad \lambda_s\ne\lambda_\ell\ {\rm{if}}\ \  r_s=r_\ell.$$  Assume without loss of generality that $r_1$ is minimal, in which case we have
 $$\sum_{\ell=1}^{m} a_\ell {\rm{ch}}_{\rm{gr}} W(\lambda_\ell,r_\ell) -u^{r_1}\sum_{\left\{\ell\ :\  r_\ell=r_1\right\}}a_\ell\ch_{\lie g}V(\lambda_\ell)\in u^{r_1+1}\bc[P][[u, u^{-1}]].$$  It follows that $$\sum_{\left\{\ell\ :\ r_\ell=r_1\right\}}a_\ell\ch_{\lie g}(V(\lambda_\ell))=0,$$ and hence, applying   Lemma \ref{ind}, we see that $a_1=0$. An obvious iteration of this argument completes the proof.

\end{pf}

\subsection{}  We need one final family of modules in $\cal I_{gr}$  and these are the  graded local Weyl modules. The module $W_{\loc}(\lambda,r)$ is the quotient of $W(\lambda,r)$  defined by imposing the additional graded  relations, $$(h\otimes t^s) w_{\lambda,r}=0,\ \ s>0.$$ The following was proved in \cite[Theorem 1]{CPweyl}, although the notation used there is slightly different. The notation used in the current paper is closer to the one in \cite{CL2}.
\begin{prop} For all $\lambda\in P^+$ and $r\in\bz$, the module $W_{\loc}(\lambda,r)$ is finite--dimensional. \hfill\qedsymbol
\end{prop}
By Lemma \ref{simpleclass} it follows that the irreducible constituents are of the form $V(\mu,s)$ for $(\mu,s)\in P^+\times\bz$ and we let  $[W_{\loc}(\lambda,r): V(\mu,s)]$ be the the multiplicity of $V(\mu,s)$ in the Jordan--Holder series in $W_{\loc}(\lambda,r)$. These multiplicities are known, see \cite{CL1} for $\lie{sl}_{r+1}$ and \cite{Naoi} for the general case and we shall give the precise statement later.

\subsection{} \label{main}  We now state our main result which connects all the modules discussed so far.  We assume that $\lie g=\lie{sl}_2$ and we  set $x_1^\pm=x^\pm$ and $h_1=h$, $\alpha_1=\alpha$ and $\omega_1=\omega$.
\begin{Theorem} Let $(m\omega,r)\in P^+\times\bz$. There exists a decreasing family $P_k\supset P_{k+1}$, $k\in\bz_+$ of objects of $\cal I_{\gr}$ such that $$P_0=P(m\omega,r),\ \ \bigcap_{k\ge 0}P_k=\{0\} $$ and  \begin{equation}\label{quot}   P_k/P_{k+1}\cong \bigoplus_{s\ge 0} W((m+2k)\omega, s )^{\oplus n(m, k,s)}\end{equation}
where $ n(m,k,s)=[W_{\loc}((m+2k)\omega,r): V(m\omega, s)].$ Moreover,
\begin{equation}\label{char}  \ch_{\gr} P(m\omega,r)=\sum_{k,s\in\bz} n(m,k,s)\ch_{\gr}W((m+2k)\omega,s).\end{equation}

\end{Theorem}

\subsection{} In general we have the following conjecture.\begin{conj}\label{bggconj}
Let $\lie g$ be a finite--dimensional simple complex Lie algebra and let $\lambda\in P^+$. There exists a decreasing family $P_k\supset P_{k+1}$, $k\in\bz_+$ of objects of $\cal I_{\gr}$ such that $$P_0=P(\lambda,r),\ \ \bigcap_{k\ge 0}P_k=\{0\} $$ and \begin{equation}\label{quotc}   P_k/P_{k+1}\cong W(\mu_k, s_k)\qquad
  \ch_{\gr} P(\lambda,r)=\sum_{k\ge 0}\ch_{\gr}W(\mu_k,s_k),\end{equation} and for all $(\mu,s)\in P^+\times\bz$ we have,  $$\#\{k: (\mu_k, s_k)= (\mu,s)\}= [W_{\loc}(\mu,r): V(\lambda, s)].$$

 \end{conj}

\subsection{}  The statement of the theorem has two parts.
The first part is that $P(m\omega,r)$ has a decreasing filtration such that \eqref{quot} is satisfied.
The second part is  that the graded  character is then given by \eqref{char} and for this, it is useful to introduce the following natural definition.

\begin{defn} Given  $M\in\Ob\mathcal I_{\rm gr}$ a Weyl flag $\cal F M$ of $M$ is a decreasing family $M_k\supset M_{k+1}$, $k\in\bz_+$ of objects of $\cal I_{\gr}$ such that $$M_0=M,\ \ \ \ \bigcap_{k\ge 0}M_k=\{0\},\ \ \  \ \ M_k/M_{k+1}\cong W(\mu_k, r_k), $$
   where $(\mu_k,r_k)\in P^+\times \bz$.  Let$$[\mathcal FM:W(\mu,r)]=\#\{k\in\bz_+: M_k / M_{k+1} \cong W(\mu, r)\}.$$
\hfill\qedsymbol\end{defn}
At this stage the only example of modules (for an arbitrary simple Lie algebra) with a Weyl flag are the global Weyl modules themselves. Notice that in this language, Theorem \ref{main} asserts that the projective modules $P(m\omega,r)$ for $\lie{sl}_2$ admit a Weyl flag.
 In Section 4 we give examples of modules admitting a Weyl flag for arbitrary simple Lie algebras. These examples also provide  some evidence for our general conjecture.

\subsection{}  The following proposition, together with Lemma \ref{linind}, shows that the multiplicity of $W(\mu,r)$ in a Weyl flag of $M\in\Ob\cal I_{\gr}$ is independent of the choice of the flag. It also   establishes \eqref{char}.

 \begin{prop}\label{weylmultdefn}  Let $\cal FM=\{M_k\}_{k\ge 0}$ be a Weyl flag of  $M\in\cal I_{\gr}$.
\begin{enumerit}
\item[(i)]  Given $s\in\bz$ there exist at most finitely many $k$ such that $M_k[s]\ne 0$.
\item [(ii)]  $$\ch_{\gr}M=\sum_{(\mu,r)\in P^+\times\bz} [\mathcal FM: W(\mu,r)]\ch_{\gr} W(\mu,r)$$
\end{enumerit}
 \end{prop}
 \begin{pf} Since $\dim M[s]<\infty$ and $M_k[s]\supseteq M_{k+1}[s]$  we see that there must exist $r$ such that $M_r[s]=M_p[s]$ for all $p\ge r$. The condition that
 $\bigcap_{k\in\bz_+}M_k=\{0\}$ now gives (i).

 For (ii), the induced filtration
 $$M[s]\supseteq M_1[s]\supseteq\cdots M_r[s]\supseteq \{0\},$$
is finite for all $s\in\bz$, and hence it follows that $$\ch_{\lie g} M[s]=\sum_{k\in\bz_+} \ch_{\lie g}W(\mu_k,s_k)[s],$$  which completes the proof of the proposition.
\end{pf}

\begin{cor} If $\cal F_1M$ and $\cal F_2M$ are two Weyl flags of $M\in {\rm Ob}\cal{I}_{gr}$ then  $$[\cal F_1M: W(\lambda,r)]= [\cal F_2M: W(\lambda,r)],\ \ (\lambda,r)\in P^+\times\bz,$$  and we denote this common multiplicity by $[M: W(\lambda,r)]$. \end{cor}
\begin{pf} This follows from part (ii) of the  proposition, together with Lemma \ref{linind}.
\end{pf}

\subsection{} As a consequence of  Lemma \ref{weyldefrel}, we see the following: for $\Gamma\subset  P^+$, we have \begin{equation}\label{weyltrunc} W(\lambda,r )^\Gamma = \begin{cases}
W(\lambda,r),\qquad  \ \lambda\in\Gamma-Q^+,\\ 0, \qquad\qquad\ \ \ \ \lambda\notin\Gamma-Q^+. \end{cases}\end{equation}
It is obvious that if $M$ has a Weyl flag and $M\in\cal I_{\gr}^\Gamma$, then $M_r$ and $M_r/M_{r+1}$ are in $\cal I^\Gamma_{\gr}$ for all $r\in\bz_+$; in other words, $M$ has  a Weyl flag in $\cal I^\Gamma_{\gr}$.
 We shall also be interested in showing that the module $P(\mu,s)^\Gamma$ admits a Weyl flag. As a first step we note the following result.

\begin{lem}\label{weylproj} Let   $\Gamma$ be a finite subset of $P^+$ and assume that $$\mu,\nu\in\Gamma\implies\ \mu\nleq\nu\ \ {\rm{and}} \ \ \nu\nleq\mu.$$ Then for all $M\in\Ob\cal I^{\Gamma}_{\gr}$ and $\nu\in\Gamma$, we have \begin{equation}\label{ninvp}\lie n^+[t]M_{\nu}=0.\  \end{equation}
 Moreover for all $(\mu,r)\in\Gamma\times \bz$, the modules    $W(\mu, r)$ are projective objects of $\cal I^{\Gamma}_{\gr}$.\end{lem}
\begin{pf}  For all  $\alpha\in R^+$ and $\mu\in\Gamma$ we have $\mu+\alpha\ge\mu$ and hence $\mu+\alpha\notin\Gamma$, which proves \eqref{ninvp}.
It follows from \eqref{weyltrunc} that for all $\mu\in\Gamma$, we have  $W(\mu,r)\in\Ob\cal I^{\Gamma}_{\gr}$.
Suppose now that we have a surjective map $\pi:M\to N$, where $M,N\in\Ob\cal I^{\Gamma}_{\gr}$ and let $\varphi: W(\mu ,r)\to N$ be   a non--zero  map of $\lie g[t]$--modules, in particular this implies that $\varphi(w_{\mu,r})\ne 0$.  Since $\pi$ is a map of $\lie g$--modules it follows that we can choose  $u\in M_{\mu}[r]$ such that $\pi(u)=\varphi(w_{\mu,r})$. The defining relations of the global Weyl modules and Equation \eqref{ninvp} imply  that the assignment $w_{\mu,r}\mapsto u$ defines a homomorphism $\psi: W(\mu,r)\to M$ of $\lie g[t]$--modules such that $\pi.\psi=\varphi$ thus proving that $W(\mu,r)$ is a projective object of $\cal I_{\gr}^\Gamma$.

\end{pf}

\section{ Characters}\label{characters}
In this section we gather some of the necessary results on global and local Weyl modules for arbitrary simple Lie algebras.  We explain what is known about the characters of these modules and state some natural identities
which are a consequence  of the  conjecture and theorem. A crucial step in understanding the modules $W(\lambda,r)$ is to determine the annihilator in $\bu(\lie h[t])$ of the element $w_{\lambda,r}$ and we deal with this in the first part of this section.
\subsection{} We  use freely the notation and definitions of the earlier sections. In addition, given $k\in\bz_+$  let $A_k$ be the polynomial ring in $k$ variables with a  $\bz$--grading given in the standard way, by setting all the generators to have grade one. The graded component  $A_k[s]$ is just the space  homogenous polynomials of degree $s$. Let  $S_k$ be  the symmetric group on $k$ letters and consider the natural graded action of   $S_k$ on $A_k$  defined by permuting the variables. Clearly $A_k[s]$ is a $S_k$--submodule of $A_k$ for all $s\in\bz$.

\subsection{} Given $\lambda\in P^+$, where $\lambda=\sum_{i=1}^nr_i\omega_i$ with $\sum_{i=1}^nr_i=r_\lambda$, we denote by $S_\lambda$ the  Young subgroup $S_{r_1}\times\cdots\times S_{r_n}$ of $S_{r_\lambda}$. The corresponding graded ring of invariants $A_{r_\lambda}^{S_\lambda}$  will be denoted as $\mathbb A_\lambda$. It is well--known that $\mathbb A_\lambda$ is again a polynomial ring in $r_\lambda$ variables with Hilbert series  $$\mathbb H(\mathbb A_\lambda)=\prod_{i=1}^n\frac{1}{(1-u)(1-u^2)\cdots (1-u^{r_i})}.$$
In  what follows, we shall regard $A_{r_\lambda}$ as the polynomial ring in the variables $t_{i,s}$, $1\le i\le n$, $1\le s\le r_i$. We shall also use the fact that $\mathbb A_\lambda$ is the subalgebra of $A_{r_\lambda}$ generated by the elements $\{(t_{i,1}^p+\cdots +t_{i,r_i}^p):\ \ 1\le i\le n,\ \ p\in\bz_+\}$.

\subsection{}  For $\lambda=\sum_{i\in I}r_i\omega_i\in P^+$  set $$\ann_{\bu(\lie h[t])} w_{\lambda,r}=\{u\in\bu(\lie h[t]): uw_{\lambda,r}=0\}.$$ It is clear that $$\ann_{\bu(\lie h[t])} w_{\lambda,r}=\ann_{\bu(\lie h[t])} w_{\lambda,0}.$$
 The following  result was established in \cite{CPweyl} in a different form, in the current language this  can be found in \cite[Lemma 8]{CFK}.

\begin{prop}\label{ann} The algebra homomorphism $\bu(\lie h[t])\to \mathbb A_\lambda$ defined by sending $$h_i\otimes t^ p\mapsto\ \ (t_{i,1}^p+\cdots +t_{i,r_i}^p),\ \ 1\le i\le n,\ \ p\in\bz,$$ induces an isomorphism of graded algebras
 $$\symm: \bu(\lie h[t])/\annh w_{\lambda,r}\cong \mathbb A_\lambda.$$ \hfill\qedsymbol
\end{prop}

\subsection{}\label{quotwel} The following is straightforward but we include a proof for the reader's convenience.
\begin{lem}Let $V\in\Ob\cal I_{\gr}$ be such that $\wt V\subset\lambda-Q^+$ and assume that $V_\lambda\ne 0$. Then $\annh w_{\lambda,r} V_\lambda=0$.  In particular,  $V_\lambda$ admits the structure of a right $\mathbb A_\lambda$--module given by $$v\boa=\boh v, \ \ v\in V_\lambda,\ \ \boh\in\bu(\lie h[t]),\ \ \boa=\symm(\boh).$$
\end{lem}
\begin{pf} Observe  that $\lie n^+[t] V_\lambda=0$ since $\wt V\subset \lambda-Q^+$. Choose  $r\in\bz$ such that $v\in V_\lambda[r]$ with $v\ne 0$.  The defining relations of $W(\lambda, r)$ given in Lemma \ref{weyldefrel} imply that there exists a non--zero morphism $W(\lambda,r)\to V$ of objects of $\cal I_{\gr}$ which sends $w_{\lambda,r}\to v$. In particular we have $uv=0$ for all $u\in\annh w_{\lambda,r}$ and the the lemma follows from Proposition \ref{ann}.
\end{pf}

\subsection{} It is  easily seen that  setting $$(yw_\lambda)u= yuw_\lambda,\ \ y\in\bu(\lie g[t]),\ \ u\in\bu(\lie h[t]),$$ defines a graded right action of $\bu(\lie h[t])$ on $W(\lambda,r)$.
 Hence we can regard $W(\lambda,r)$ as a graded  right module for $\mathbb A_\lambda$ and in fact as a graded $(\lie g[t], \mathbb A_\lambda)$--bimodule. In this language one has the following equivalent definition \cite{CFK} (see also \cite[Section 4.3] {BCGM}) of the $W_{\loc}(\lambda,r)$, namely:
\begin{equation}\label{altlocal} W_{\loc}(\lambda,r)\cong_{\lie g[t]} W(\lambda,r)\otimes_{\mathbb A_\lambda}\bc_\lambda,
\end{equation}
where $\bc_\lambda$ is the one dimensional module obtained by going modulo the augmentation ideal in $\mathbb A_\lambda$. The following is immediate.
 \begin{lem}\label{invn}  Given $\lambda,\mu\in P^+$, the subspace $W(\lambda,r)_\mu^{\lie n^+}$ is a graded $\mathbb A_\lambda$--submodule of $W(\lambda,r)$ and we have an isomorphism of $(\lie g,\mathbb A_\lambda)$--bimodules \begin{gather*}\bu(\lie g)W(\lambda,r)^{\lie n^+}_\mu\cong V(\mu)\otimes_\bc W(\lambda,r)^{\lie n^+}_\mu, \\  \bu(\lie g)W_{\loc}(\lambda,r)^{\lie n^+}_\mu\cong V(\mu)\otimes_\bc( W(\lambda,r)^{\lie n^+}_\mu\otimes_{\mathbb A_\lambda} \bc_\lambda),\end{gather*} which sends $$w\to v_\mu\otimes w,\ \ \ (resp.\   w\to v_\mu\otimes(w\otimes 1)), \  \ w\in W(\lambda,r)^{\lie n^+}_\mu.$$  Here we regard $\lie g$ as acting only on $V(\mu)$ and $\mathbb A_\lambda$ on $W(\lambda,r)^{\lie n^+}_\mu$.

\hfill\qedsymbol
 \end{lem}
  Since $W(\lambda,r)$ is generated as a $\lie g$--module by the spaces $W(\lambda,r)_\mu^{\lie n^+}$, we have an isomorphism of $(\lie g,\mathbb A_\lambda)$--bimodules, \begin{equation}\label{bimod}W(\lambda,r)\cong \bigoplus_{\mu\in P^+}V(\mu)\otimes W(\lambda,r)^{\lie n^+}_\mu,\end{equation} and hence we see that
$$\ch_{\gr}W(\lambda,r)=\sum_{\mu\in P^+}\mathbb H(W(\lambda,r)^{\lie n^+}_\mu)\ch_{\lie g} V(\mu).$$

 \subsection{} To determine the graded character of $W(\lambda,r)$ we need the following result which was established in \cite{CPweyl} for $\lie{sl}_2$, in \cite{CL1} for $\lie{sl}_{r+1}$, in \cite{FoL} for simply--laced $\lie g$ and in \cite{Naoi} in general. It can also be deduced from the quantum case through the work of  \cite{BN}, \cite{Ka} and \cite{Nak}.
\begin{thm} Regarded as a module for $\mathbb A_\lambda$, the global Weyl module $W(\lambda,r)$ is free of rank $\dim_{\bc} W_{\loc}(\lambda,r)$.\hfill\qedsymbol
\end{thm}
We note the following corollary.
\begin{cor}\label{nfree} For $\mu\in P^+$, the subspace $W(\lambda,r)_\mu^{\lie n^+}$ is a free $\mathbb A_\lambda$--module of rank $\dim_{\bc}W_{\loc}(\lambda,r)_\mu^{\lie n^+}$.\end{cor}
\begin{pf} Using \eqref{ndecomp} and  \eqref{bimod} we see that $W(\lambda,r)_\mu^{\lie n^+}$ is an $\mathbb A_\lambda$ summand of $W(\lambda,r)$ and hence is a projective $\mathbb A_\lambda$--module. Since $\mathbb A_\lambda$ is a polynomial ring it follows from the the famous result of Quillen (see \cite{Lam} for an exposition) that $W(\lambda,r)^{\lie n^+}_\mu$ is in fact a free $\mathbb A_\lambda$--submodule of $W(\lambda,r)$.

 To determine its rank, let  $\pi_\lambda: W(\lambda,r)\to W_{\loc}(\lambda,r)$ be the canonical morphism of graded $\lie g[t]$--modules sending $w_{\lambda,r}\to w_{\lambda,r}\otimes 1$ where we use the equivalent  definition given in \eqref{altlocal}. Then, $\pi_\lambda(W(\lambda,r)_\mu^{\lie n^+})= W_{\loc}(\lambda,r)_\mu^{\lie n^+}$ and hence $$\rk_{\mathbb A_\lambda}W(\lambda,r)_\mu^{\lie n^+}\ =\dim_{\bc} W_{\loc}(\lambda,r)_\mu^{\lie n^+}.$$
\end{pf}

\subsection{}\label{grw} The following proposition determines  the  graded character of $W(\lambda,r)$ in terms of the local Weyl module $W_{\loc}(\lambda,r)$.

\begin{prop}\label{grchar} For $(\lambda,r)\in P^+\times\bz$, we have $$\ch_{\gr}W(\lambda,r)=\ch_{\gr}W_{\loc}(\lambda,r)\mathbb H(\mathbb A_\lambda).$$

\end{prop}
\begin{pf} Using  \cite[Proposition 11.4.7]{Bo} we see that Corollary \ref{nfree} implies that the $\mathbb A_\lambda$--module $W(\lambda,r)^{\lie n^+}_\mu$ is graded free. In other words, if we  pick a graded basis of $W_{\loc}(\lambda,r)^{\lie n^+}_\mu$, then   a set of graded pre-images  is a graded basis for   $W(\lambda,r)^{\lie n^+}_\mu$ as an $\mathbb A_\lambda$--module.
This implies that\begin{gather*}\sum_{s\in\bz}\dim W(\lambda,r)^{\lie n^+}_\mu[s]u^s =\left(\sum_{s\in\bz}\dim W_{\loc}(\lambda,r)^{\lie n^+}_\mu[s]u^s\right)\mathbb H(\mathbb A_\lambda),\end{gather*} and
hence we get \begin{gather*}\ch_{\gr}W(\lambda,r) =\sum_{\mu\in P^+}\ch_{\lie g}V(\mu)\left(\sum_{s\in\bz}\dim W_{\loc}(\lambda,r)^{\lie n^+}_\mu[s]u^s\right)\mathbb H(\mathbb A_\lambda) =\ch_{\gr}W_{\loc}(\lambda,r)\mathbb H(\mathbb A_\lambda).\end{gather*}
\end{pf}

\subsection{}\label{grloc} The graded character of the local Weyl module was computed in \cite{CL1} for $\lie{sl}_{r+1}$ and in \cite{Naoi} in general. We shall just give the result for $\lie{sl}_{r+1}$ since the general case requires the introduction of a number of ideas from the theory of crystal bases.

Regard  $\lambda=\sum_{i=1}^nr_i\omega_i\in P^+$  as an $n$--tuple of integers $\lambda=(r_1,\cdots, r_n)$.
Let $\xi=(\xi_1\ge\cdots\ge\xi_{n+1}\ge 0)$ be  a partition of $r_\lambda=\sum_{i=1}^{n}r_i$ with at most $(n+1)$--parts and let  $K_{\lambda,\xi^{tr}}(u)$ be  the generalized Kostka polynomial defined in \cite[III.6]{Ma}. Set $\mu_\xi=\sum_{i=1}^n(\xi_i-\xi_{i+1})\omega_i.$

\begin{prop}\label{chgrloc} For $\lambda\in P^+$ we have,$$\ch_{\gr} W_{\loc}(\lambda,r)=u^r\sum_\xi K_{\lambda,\xi^{tr}}(u)\ch V(\mu_\xi),$$ where the sum is over all partitions $\xi$ of $r_\lambda$ with at most $(n+1)$--parts.
\hfill\qedsymbol
\end{prop}

\subsection{} Assuming that Conjecture \ref{bggconj} is true, we now compute the character of the projective module $P(\lambda,0)$  in terms of the local Weyl module as follows.     Using Proposition \ref{grw} we have  \begin{eqnarray*}&\ch_{\gr}P(\lambda,0)&=\sum_{(\mu,s)\in P^+\times\bz}[P(\lambda,0):W(\mu,s)]\ch_{\gr} W(\mu,s)\\ \\ &&=\sum_{(\mu,s)\in P^+\times\bz}[W_{\loc}(\mu,0): V(\lambda,s)]u^s\ch_{\gr} W(\mu,0)\\ &&=\sum_{\mu\in P^+}\ch_{\gr} W(\mu,0)\sum_{s\in\bz}K_{\mu,\xi_\lambda^{tr}}(u)_su^s\\ &&=\sum_{\mu\in P^+}\ch_{\gr} W(\mu,0)K_{\mu,\xi_\lambda^{tr}}(u) \\
&&= \sum_{\mu\in P^+(\lambda)}K_{\mu,\xi_\lambda^{tr}}(u)\ch_{\gr} W_{\loc}(\mu,0)\mathbb H(\mathbb A_\mu),\end{eqnarray*}
where $P^+(\lambda)$ consists of $\mu\in P^+$ such that there exists $\xi_\lambda^{tr}=(\xi_1\ge\cdots\ge \xi_{n+1})$  a partition of $r_\mu$ into at most $(n+1)$--parts and $\lambda=\sum_{i=1}^n(\xi_i-\xi_{i+1})\omega_i$.

\subsection{} A simple consequence of the Poincare--Birkhoff--Witt theorem (see \cite[Proposition 2.1]{CG1}) is that we have an isomorphism of $\bz$--graded $\lie g[t]$--modules,  $$P(\lambda,r)\cong_{\lie g[t]}\bu(\lie g\otimes t\bc[t])\otimes_\bc V(\lambda,r).$$ Here we regard $\bu(\lie g\otimes t\bc[t])\otimes V(\lambda, r)$ as a $\lie g[t]$--module as follows: 
$$(x\otimes t^s)(y\otimes v)=\begin{cases}(x\otimes t^s)y\otimes v,\ \ s>0,\\ [x,y]\otimes v+ y\otimes xv, \ \ \ s=0,\end{cases}$$
 where $x\in\lie g$, $y\in\bu(\lie g\otimes t\bc[t])$, $s\in\bz_+$ and $v\in V(\lambda,r)$. In particular $$P(0,r)\cong_{\lie g[t]} \bu(\lie g\otimes t\bc[t])\otimes V(0,r),\ \ \ \ch_{\gr}P(\lambda,r)=u^r\ch_{\lie g}V(\lambda)\ch_{\gr}P(0,0),$$
and we get
$$\sum_{\mu\in P^+(\lambda)}u^rK_{\mu,\xi_\lambda^{tr} }(u)\ch_{\gr} W_{\loc}(\mu,0)\mathbb H(\mathbb A_\mu)=\sum_{\mu\in P^+(0)}u^rK_{\mu,\xi_0^{tr} }(u)\ch_{\gr} W_{\loc}(\mu,0) \ch_{\lie g} V(\lambda)\mathbb H(\mathbb A_\mu).$$

\subsection{}\label{kostkasl2}  By the Poincare--Birkhoff--Witt theorem we  have an isomorphism of graded vector spaces  $$\bu(\lie g\otimes t\bc[t])\cong S(\lie g\otimes t\bc[t]),$$ where for a vector space $V$ we denote by $S(V)$ the symmetric algebra of $V$. Consider now the special case when $\lambda=0$ and $\lie g$ is $\lie{sl}_2$.
The preceding discussion, together with Theorem 1, proves the following:
\begin{Theorem} Suppose that $\lie g=\lie{sl}_2$.
Then, $$\ch_{\gr}S(\lie g\otimes t\bc[t])=\sum_{r,m}K_{2m, (m\ge m) }(u)K_{2m, (2m-r\ge r)}(u)\mathbb H(\mathbb A_{2m})\ch_{\lie g}V(2m-2r).$$ Considering the dimensions of the graded pieces on both sides, we have
$$\prod_{r\ge 1}\left(\frac{ 1}{1-u^r}\right)^{\dim\lie g}=\sum_{r,m}(2m-2r+1)K_{2m, (m\ge m) }(u)K_{2m, (2m-r\ge r)}(u)\mathbb H(\mathbb A_{2m}).$$
 \hfill\qedsymbol
\end{Theorem}
\section{Examples of Modules admitting a Weyl flag}\label{examples}
In this section we give a  non--trivial example  of an object of $\cal I_{\gr}$ which admits a Weyl flag. We begin with the following remark. In the case of the BGG category $\cal O$ for $\lie g$, it is relatively easy to produce examples of modules admitting a Verma flag. The standard example of such a module is the tensor product of a Verma module with a finite--dimensional module for $\lie g$. In the case of $\cal I_{\gr}$ such a statement is false. As an example, suppose that $\lie g$ is $\lie{sl}_2$ and consider the tensor product  $V=W(\omega,0)\otimes V(\omega,0)$. Suppose that $V$ admits a Weyl flag.  Since $V_{2\omega}[0]\ne 0$ and $\wt V\subset\{\pm 2\omega,0\},$ we see that $W(2\omega,0)$ must be a sub-quotient of $V$. Using  the results of Section \ref{characters}, we get $$\mathbb H(W(2\omega,0)_{2\omega}) =\mathbb H(\mathbb A_{2\omega}).$$ On the other hand $$\mathbb H(V_{2\omega})= \mathbb H\left( (W(\omega,0)\otimes V(\omega,0))_{2\omega}\right)=\mathbb H(\mathbb A_\omega).$$ Since $\mathbb A_{2\omega}\cong A_2^{S_2}$ and $\mathbb A_\omega=A_1$ we get $\dim W(2\omega,0)_{2\omega}[2]=2> \dim V_{2\omega}[2].$ This is clearly impossible and hence $V$ does not admit a Weyl flag.

\subsection{} Let $\theta\in R^+$ be the highest root of $\lie g$ and recall that $\theta\in P^+$. We shall prove,
\begin{prop}\label{firststep} For $r\in\bz$, we have $$\Ext^1_{\cal I_{\gr}}(\tau_r\bc, W(\theta,r+1))\ne 0.$$ \end{prop}
\begin{cor}
 We have a  non--split short exact sequence of objects $$0\to W(\theta, r+1)\to P(0,r)^\theta\to\tau_r\bc\to 0$$ in the category $\mathcal I_{\gr}^\theta$.
 \end{cor}
The proofs of the proposition and its corollary are given in the rest of this section.

 \subsection{}\label{weyltheta} The structure of the local Weyl module $W_{\loc}(\theta, 0)$ can be read off from the character formulae discussed in the previous sections. In this special case, they can also be found in \cite{CM}  in the following explicit way. We have $V(\theta)\cong\lie g$ where we regard $\lie g$ as a module for itself through  the adjoint representation.  Let $\langle\cdot , \cdot\rangle $ be the Killing form on $\lie g$.
Then, $$W_{\loc}(\theta,0)[s]=0,\  s\ge 2,\quad \ W_{\loc}(\theta,0)[0]\cong_{\lie g} \lie g,\ \ \quad W_{\loc}(\theta,0)[1]\cong_{\lie g[t]} \tau_1\bc,$$ and the action of $\lie g\otimes \bc[t]$ is $$(x\otimes t^s)(y,a)= 0,\ \ s\ge 2,\ \ (x\otimes t)(y,a)= (0,\langle x,y\rangle),\ \ \ x(y,a)=([x,y],0).$$
We remark for the reader's convenience that the element $w_{\theta,0}$ is just $x^+_\theta\in \lie g$.
As a consequence of Proposition \ref{grchar}, we get \begin{equation}\label{grchartheta}\ch_{\gr}W(\theta,0)=\ch_{\gr}W_{\loc}(\theta,0)\mathbb H(\mathbb A_\theta)=(\ch_{\lie g}V(\theta)+u\ch_{\lie g}\bc)\mathbb H(\mathbb A_\theta).\end{equation}

\subsection{} Define the module $\widehat {W}(\theta,0)\in\cal I_{\gr}$ as follows: $$\widehat{W}(\theta,0)= W_{\loc}(\theta,0)\otimes\bc[t]= (\lie g\oplus \tau_1\bc)\otimes\bc[t],$$
with the $\lie g[t]$ action given by\begin{equation}\label{hatweq} (x\otimes t^r)(y\otimes f,a\otimes g)= ([x,y]\otimes t^rf, 0) + r(0,\langle x ,y\rangle\otimes t^{r-1}f)\end{equation} where $x,y\in\lie g$,\ \ $a\in\bc$, and $f,g\in \bc[t]$. Clearly,\begin{equation}\label{grhat}\ch_{\gr}\widehat{W}(\theta,0)=\ch_{\gr}W_{\loc}(\theta,0)\mathbb H(A_1).\end{equation}
\begin{lem}\label{hatw} We have $\widehat{W}(\theta,0)\in\Ob\cal I_{\gr}^\theta$. Moreover it is generated by the element $x^+_\theta\otimes 1$ and hence is a quotient of $W(\theta,0)$. Furthermore, if $\lie g$ is not of type $A_n$ or $C_n$, then $$\widehat{W}(\theta,0)\cong W(\theta,0).$$
\end{lem}
\begin{pf} The first statement is clear since $$\wt(\widehat{W}(\theta,0))=\wt(W_{\loc}(\theta,0))=\wt V(\theta)\subset \theta-Q^+.$$ Hence $\lie n^+[t]\ x^+_\theta =0,$ and the second statement follows from the defining relations of $W(\theta,0)$ once we prove that $x_\theta^+\otimes 1$ generates $\widehat W(\theta,0)$ as  $\lie g[t]$--module. But this is easily checked, given the explicit formulae in \eqref{hatweq}.

Suppose now that $\lie g$ is not of type $A_n$ or $C_n$. Then $\theta\in\{\omega_i: i\in I\}$ and hence $\mathbb A_\theta\cong A_1$ as graded algebras. Equations \eqref{grchartheta}, \eqref{grhat} imply that $\widehat{W}(\theta,0)$ and $W(\theta,0)$ have the same graded characters and the final statement of the lemma is now immediate.
\end{pf}
 \subsection{} \label{nonzeroext} The following lemma, together with Lemma \ref{hatw}, completes the proof of Proposition \ref{firststep} in the case when $\lie g$ is not of type $A_n$ or $C_n$.
\begin{lem} For $r\in\bz$, we have $\Ext^1_{\cal I_{\gr}^\theta}\left(\tau_r\bc, \ \tau_{r+1}\widehat{W}(\theta,0)\right)\ne 0$.
\end{lem}
\begin{pf} Define an action of $\lie g[t]$ on the direct sum of $\lie g$--modules $$U=\tau_r\bc\oplus \tau_{r+1}\widehat{W}(\theta,0),$$ as follows: $\tau_{r+1}\widehat{W}(\theta,0)$ is a $\lie g[t]$--submodule of $U$, and the action of $\lie g[t]$ on $\tau_r\bc$ is given by $$(x\otimes t^s)(a,0,0)=as(0, x\otimes t^{s-1}, 0),\ \ s\in \bn,\ \  x\in\lie g.$$ A simple checking shows that this is indeed an action of $\lie g[t]$ and also that $U$ is  indecomposable.  In other words, we have a
non--split short exact sequence $$0\to \tau_{r+1}\widehat W(\theta,0)\to U\to \tau_r\bc\to 0,$$   and the lemma is proved.
\end{pf}

\subsection{} \label{K}
Let $K$ be the kernel of the map $P(0,r)^\theta\to \tau_r\bc\to 0$. To complete the proof of Proposition \ref{firststep} and its corollary, we will prove that $$K\cong W(\theta,r+1),$$ as objects of $\cal I_{\gr}$.   We begin with the following.

\begin{lem}   We have $$K=\bu(\lie g[t])(x^+_\theta\otimes t)p_r,\qquad {\rm{where}}\ \ p_r=p_{0,r}.$$ In particular, there exists a surjective morphism $\psi_r: W(\theta, r+1)\to K\to 0$ of objects of $\cal I_{\gr}$ with $\psi_r(w_{\theta,r+1})=(x^+_\theta\otimes t)p_r.$
\end{lem}
\begin{pf} Since $\lie g \ p_r=0$ and $P(0,r)^\theta[r]=\bc p_r$ (see Lemma \ref{pdefrel})
 it follows  that $$P(0,r)^\theta=\bu(\lie g\otimes t\bc[t])p_r,\qquad \ K=\bu(\lie g\otimes t\bc[t])_+p_r. $$  The elements $\{x\otimes t: x\in\lie g\}$ generate
 $\lie g\otimes t\bc[t]$ as a Lie algebra and hence  generate $\bu(\lie g\otimes t\bc[t])$ as an algebra. Hence it suffices to show that $$(x\otimes t)p_r\in\bu(\lie g[t])(x^+_\theta\otimes t)p_r,\ \ \ {\rm{for\ all}}\ \ x\in\lie g.$$ Since $$y(x\otimes t)p_r= ([y,x]\otimes t)p_r,\ \ \ x,y\in\lie g,$$ it is enough to prove that the $\lie g$--module (via the adjoint action of $\lie g$)  through $(x^+_\theta\otimes t)$ is
all of $\lie g\otimes t$. But this is equivalent to the well--known statement \cite{Hu} that $x^+_\theta$ generates the adjoint representation of $\lie g$. Since $\wt K\subset \theta-Q^+$, we see that  $\lie n^+[t]((x^+_\theta\otimes t)p_r)=0$ and hence $K$ is a quotient of $W(\theta,r+1)$.

\end{pf}

  \subsection{} Suppose that $\lie g$ is not of type $A_n$ or $C_n$, in which case we only have to prove Corollary \ref{firststep}. Consider the short exact sequence $$0\to  K\to P(0,r)^\theta\to\tau_r \bc\to 0.$$ Since $W(\theta, r+1)[r]=0$ we have,  $$\Hom_{\cal I^\theta_{\gr}}(\tau_r\bc, W(\theta,r+1))=\Hom_{\cal I^\theta_{\gr}}(P(0,r)^\theta, W(\theta, r+1))=0,  $$
and since  $P(0,r)^{\theta}$ is projective in $\cal I^\theta_{\gr}$ we get by using Lemma \ref{nonzeroext} that $$\Hom_{\cal I^\theta_{\gr}}(K, W(\theta,r+1))\cong \Ext^1_{\cal I^{\theta}_{\gr}}(\tau_r\bc, W(\theta, r+1))\ne 0.$$
  Let $\varphi\in\Hom_{\lie g[t]}( K, W(\theta,r+1))$ be   non-zero.
    By Lemma \ref{K}, the element $(x_\theta^+ \otimes t)p_r$ spans $K_\theta[r+1]$ and hence we may assume  $$\varphi((x_\theta^+\otimes t)p_r)= w_{\theta,r+1}.$$ But now, $\psi_r.\varphi:W(\theta,r+1)\to W(\theta,r+1)$ is a surjective morphism which splits since $W(\theta,r+1)$ is projective (see Lemma \ref{weylproj}) in $\cal I^\theta_{\gr}$. Since $W(\theta,r+1)$ is indecomposable it follows that $\varphi$ is an isomorphism and Corollary \ref{firststep} is proved in this case.

\subsection{} To deal with the case when $\lie g$ is of type $A_n$ or $C_n$, we need the following result which was proved in \cite[Theorem 3]{BCGM}.  In that paper it was shown that the space of all $\lie g[t]$--maps (not just the graded ones) between global Weyl modules $W(\lambda,r)$ and $W(\mu,s)$ is zero for $\lambda\ne \mu$ when $\lie g$ is of type $A_n$ or $C_n$.
\begin{lem}\label{homs} Assume that $\lie g$ is of type $A_n$ or $C_n$. Then, for $r,s\in\bz$ we have $$\Hom_{\cal I_{\gr}}(\tau_r\bc, W(\theta,s))=0.$$ Further, any non--zero element of $\Hom_{\cal I_{\gr}}(W(\theta, r), W(\theta, s))$ is injective.

\hfill\qedsymbol
\end{lem}
\subsection{} Suppose now that $\lie g$ is of type $A_n$ or $C_n$, in which case  we shall prove Proposition \ref{firststep} and its corollary simultaneously.
Consider the induced sequence 
$$0\to\Hom_{\cal I^\theta_{\gr}}(\tau_r\bc, W(\theta,r-1))\to \Hom_{\cal I^\theta_{\gr}}(\P(0,r)^\theta, W(\theta,r-1))\to\Hom_{\cal I^\theta_{\gr}}(K, W(\theta,r-1)).$$ 
Lemma \ref{homs} and Equation \ref{grchartheta} give 
$$\Hom_{\cal I^\theta_{\gr}}(\tau_r\bc, W(\theta,r-1))=0,\qquad [W(\theta, r-1)[r]: \tau_r\bc]=1,$$ 
and hence we have $W(\theta,r-1)[r]^{\lie n^+}_0\ne 0$. The defining relations of $P(0,r)$  show that 
$$\Hom_{\cal I^\theta_{\gr}}(P(0,r)^\theta, W(\theta,r-1))\ne 0$$
which forces 
$$\Hom_{\cal I^\theta_{\gr}}( K, W(\theta,r-1))\ne 0.$$ 
Let $\varphi\in \Hom_{\cal I^\theta_{\gr}}( K, W(\theta,r-1))$ be   non--zero, in which case $\varphi((x^+_\theta\otimes t)p_r)\ne 0$. Then the composite map $\varphi.\psi_r: W(\theta,r+1)\to W(\theta,r-1)$ is non--zero and hence injective by Lemma \ref{homs}. It follows that $\varphi$ is injective, and hence, again using Lemma \ref{homs}, that $$K\cong\Im\varphi\cong W(\theta,r+1).$$

\section{Proof of Theorem 1: the first reduction}\label{red1} The proof of Theorem 1 requires an explicit realization of the global Weyl modules and we begin this section by recalling this construction. We note here that it suffices, by application of the functor $\tau_r$, to prove Theorem 1 in the case of $P(\ell\omega,0)$.

\subsection{} Throughout the next two sections we will be working with $\lie{sl}_2$ and we recall that we write $x^\pm$ for $x_1^\pm$ etc. We shall use without further mention the fact that in this case  the Lie subalgebras $\lie h[t]$ and $\lie n^\pm[t]$  are abelian. The corresponding universal enveloping algebras are polynomial rings in the  infinitely many variables $\{x^\pm\otimes t^s: s\in\bz_+\}$. For  $s\in\bz_+$, and  $\bor=(r_1,\cdots, r_s)\in\bz_+^s$,
set $$\bx^\pm_\bor=(x^\pm\otimes t^{r_1})\cdots(x^\pm\otimes t^{r_{s}})$$ The elements $\{\bx^+_\bor:\bor\in\bz_+^s\}$ are a basis for $\bu(\lie n^+[t])_{s\alpha}$, where we understand $\bx^\pm_\bor=1$ for $\bor\in\bz^0_+$.

\subsection{}
 Let   $e_\pm$ be the standard basis for $V(\omega)$ satisfying $$x^\pm e_\pm=0,\  \ \  x^\pm e_\mp=e_\pm,\ \ \ he_\pm=\pm e_\pm.$$
Given $\ell \in\bz_+$, define a structure of a $(\lie{sl}_2[t], A_\ell)$-bimodule on the space
 $V(\omega)^{\otimes \ell}\otimes A_\ell$, by extending linearly the assignment, \begin{gather*} (z\otimes t^r).
  (v_1\otimes\cdots \otimes v_\ell\otimes f)=\sum_{j=1}^\ell v_1\otimes\cdots\otimes v_{j-1}\otimes z.v_j\otimes v_{j+1}\otimes\cdots\otimes v_\ell\otimes t_j^rf,\\
 (v_1\otimes\cdots\otimes v_\ell\otimes f).g=(v_1\otimes\cdots\otimes v_\ell\otimes fg),\\ \end{gather*}where $ z\in\lie g$, $v_s\in V(\omega)$, $1\le s\le \ell$, $f,g\in A_\ell$ and $ r\ge 0.$

 The space $(V(\omega)^{\otimes \ell}\otimes A_\ell)$ inherits  the grading on $A_\ell$, i.e., for an integer $r$, the $r^{th}$--graded piece is $(V(\omega)^{\otimes \ell})\otimes A_\ell[r]$.
The following is trivial.
\begin{lem} For all $\ell\in\bz_+$ we have $V(\omega)^{\otimes \ell}\otimes A_\ell\in\Ob\cal I_{gr}^{\ell\omega}$ .\hfill\qedsymbol\end{lem}

\subsection{}  The symmetric group $S_\ell$ acts  on $V(\omega)^{\otimes \ell}$ by permuting the factors in the tensor product. The induced diagonal action of $S_\ell$ on $V(\omega)^{\otimes \ell}\otimes A_\ell$
commutes with  the right action of $A_\ell$ and $\lie{sl}_2[t]$. In particular, this means that $(V(\omega)^{\otimes \ell}\otimes A_\ell)^{S_\ell}\in\Ob\cal I_{\gr}$. Recall from Section 3 that  $$\mathbb A_{\ell\omega}= A_\ell^{S_\ell},$$  and hence $W(\ell\omega,r)$ is a right module for $A_\ell^{S_\ell}$.
The following was proved in \cite{CPweyl} in the case that $\boa=1$; the proof in the general case is identical.
\begin{thm}\label{globalweyl}The assignment $w_{\ell,0}\to e_+^{\otimes \ell}\otimes \boa$, $\boa\in A_\ell^{S_\ell}$, induces an isomorphism of graded $(\lie g[t], A_\ell^{S_\ell}$)--bimodules, $$W(\ell\omega, \gr \boa)\cong (V(\omega)^{\otimes \ell}\otimes \boa A_\ell)^{S_\ell}.$$

 \hfill\qedsymbol
\end{thm}

\subsection{} A much more general version of the following result was established in \cite[Proposition 6.1]{BCGM}.
\begin{prop} Let $\bov\in (V(\omega)^{\otimes \ell}\otimes A_\ell)$. Then  $$(x^+\otimes\bc[t])\bov=0 \ \ \ \iff\ \ \bov\in (e_+^{\otimes \ell}\otimes A_\ell).$$
\end{prop}
Together with Theorem \ref{globalweyl} we have the following,
\begin{cor}\label{hw} Let $w\in W(\ell\omega, 0)$. Then,$$(x^+\otimes\bc[t])w=0 \ \ \iff \ \  w\in W(\ell\omega, 0)_{\ell\omega}.$$\end{cor}
\subsection{} Set $e_+\wedge e_-=(e_+\otimes e_- -e_-\otimes e_+)$ and for $k\in\bz_+$, denote by  $( e_+\wedge e_-)^{\otimes k}$ the element $(e_+\otimes e_- -e_-\otimes e_+)^{\otimes k}$ of $V(\omega)^{\otimes 2k}$.

 \begin{prop} For  $m,k\in\bz_+$,
the assignment $$p_{m,0}\mapsto e_+^{\otimes m}\otimes (e_+\wedge e_-)^{\otimes k}\otimes 1$$ defines a homomorphism  $ \psi_{m,k}: P(m\omega,0)\to V(\omega)^{\otimes(m+2k)}\otimes A_{m+2k}$ of objects of $\cal I_{\gr}$.
\end{prop}
 \begin{pf}
It is simple to check that \begin{gather*}h\left(e_+^{\otimes m}\otimes (e_+\wedge e_-)^{\otimes k}\right) = m(e_+^{\otimes m}\otimes (e_+\wedge e_-)^{\otimes k}),\\ \ x^+\left(e_+^{\otimes m}\otimes (e_+\wedge e_-)^{\otimes k}\right)=0.\end{gather*} The proposition is now a consequence of the defining relations of $P(m\omega,0)$.\end{pf}

\begin{cor}
 The map $\psi_{m,k}$ induces a map $\psi^{k}_m: P(m\omega,0)^{(m+2k)\omega}\to V(\omega)^{\otimes (m+2k)}\otimes A_{m+2k}$ of objects of $\cal I^{(m+2k)\omega}_{\gr}$.\hfill\qedsymbol\end{cor}

\subsection{} {\em {For the rest of the paper we fix $m,k\in\bz_+$ and set $\ell=m+2k$. To simplify notation, we shall suppress $\omega$; for example, we write $P(m, 0)^{\ell}$ for $P(m\omega,0)^{\ell\omega}$. We shall also write $p_m$ for $p_{m,0}$}}.

\vskip24pt

\begin{lem}\label{genpl} The subspace $P(m,0)^\ell_\ell$ is a left  $\lie h[t]$--submodule of $P(m,0)^\ell$ and is  generated by elements of the set $\{\bx^+_\bor p_m:\bor\in\bz_+^k\}$.
\end{lem}
\begin{pf} Since $\lie h\otimes\bc[t]$ is abelian it is clear that $P(m,0)^\ell_\ell$ is a $\lie h[t]$--submodule of $P(m,0)^\ell$. It is also trivial to see that the elements $\bx^+_\bor p_m\in P(m,0)^\ell_\ell$. Since $$u\in\left(\lie n^-\bu(\lie g[t])\right)_{k\alpha}\implies u\in\sum_{r>k}\bu(\lie g[t])(\bu(\lie n^+[t]))_{r\alpha},$$ and $\wt P(m,0)^\ell\subset \ell\omega-Q^+$, we see that $$u\in\left(\lie n^-\bu(\lie g[t])\right)_{k\alpha}\implies up_m=0,$$
 and hence $$P(m,0)^\ell_\ell=\sum_{\bor\in\bz_+^k}\bu(\lie h[t])\bx_\bor p_m,$$ thus  completing the proof of the lemma.
\end{pf}

\subsection{}

For $\bor\in\bz_+^k$, define elements $\bop(\bor)\in A_\ell$ by,
\begin{equation}\label{defnp}
\bop(\bor)=\sum_{\sigma\in
S_{k}}(t_{m+1}^{r_{\sigma(1)}}-t_{m+2}^{r_{\sigma(1)}})\cdots(t_{m+2k-1}^{r_{\sigma(k)}}-t_{m+2k}^{r_{\sigma(k)}}).
\end{equation}

\begin{lem} \label{impsi} For $\bor\in\bz_+^k$ and $s\in\bz_+$, we have $$\psi^k_m(\bx_\bor p_m)=e_+^{\otimes \ell}\otimes\bop(\bor),\qquad  \psi^k_m((h\otimes t^s)\bx_\bor p_m)=e_+^{\otimes \ell}\otimes\bop(\bor)(t_1^s+\cdots+t_\ell^s).$$
 \end{lem}
\begin{pf} The lemma is proved by a straightforward calculation. We illustrate this in the case $k=1$. Thus for $r\in\bz_+$, we have \begin{eqnarray*}&(x\otimes t^r)(e_+^{\otimes m}\otimes (e_+\otimes e_- - e_-\otimes e_+)\otimes 1)&= e_+^{\otimes m}(e_+\otimes e_+\otimes t^r_{m+1}-e_+\otimes e_+\otimes t_{m+2}^r)\\
&&= e_+^{\otimes (m+2)}\otimes (t^r_{m+1}-t^r_{m+2}),\end{eqnarray*} as needed.

\end{pf}

\subsection{} Let $M_{k,\ell}$ be the $A_\ell^{S_\ell}$--submodule of $A_\ell$ generated by the elements $\{\bop(\bor):\bor\in\bz_+^k\}$. Then Lemma \ref{genpl} and Lemma \ref{impsi} together yield,
\begin{equation}\label{impsi2} \psi^k_m(P(m,0)^\ell_\ell)= e_+^{\otimes\ell}\otimes M_{k,\ell}.\end{equation}
The module $M_{k,\ell}$ was studied in \cite[Theorem 1] {CL2} where the following result was established.

\begin{thm}\label{clfree} We have an isomorphism of graded $A_{\ell}^{S_\ell}$--modules $$W(\ell,0)^{\lie n^+}_{m}\cong M_{k,\ell},$$ where we recall that $m=\ell-2k$. In particular, $M_{k,\ell}$ has an $A_\ell^{S_\ell}$ basis consisting  of graded elements $\{\boa_1,\cdots ,\boa_N\}$ where $N=\dim_{\bc}W_{\loc}(\ell,0)_m^{\lie n^+}$ and the number of basis elements of a fixed grade $s$ is $[W_{\rm loc}(\ell,0):V(m,s)]$.
\hfill\qedsymbol
\end{thm}

\subsection{} Assume now that we have fixed a basis $\{\boa_1,\cdots,\boa_N\}$ of $M_{k,\ell}$ as in Theorem \ref{clfree} and let $\gr\boa_r$ be the grade of $\boa_r$, $1\le r\le N$.

\begin{prop}\label{direct} The $\lie g[t]$--submodule $\bw(\ell)$ of $V(\omega)^{\otimes \ell}\otimes A_\ell$ which is generated by elements of  the set $\{e_+^{\otimes\ell}\otimes \bop(\bor): \bor\in\bz_+^k\}$ is an object of $\cal I_{\gr}$ and $$\bw(\ell)\cong \bigoplus_{r=1}^NW(\ell , \gr\boa_r).$$\end{prop}
\begin{pf} The elements $\bop(\bor)$ are graded elements of $A_\ell$ which shows that the $\lie g[t]$--submodule generated by the elements $\{e_+^{\otimes\ell}\otimes \bop(\bor): \bor\in\bz_+^k\}$ is a graded submodule  of $V(\omega)^{\otimes \ell}\otimes A_\ell$ and hence an object of $\cal I_{\gr}$. For $1\le r\le N$, the elements $e_+^{\otimes \ell}\otimes \boa_r$ satisfy the defining relations of $w_{\ell, \gr\boa_r}$. Since these elements have the same grade, this proves that there exists a surjective map of graded $\lie g[t]$--modules $$\pi: \bigoplus_{r=1}^NW(\ell ,\gr\boa_r)\twoheadrightarrow \bw(\ell),\qquad \ w_{\ell,\gr\boa_r}\to e_+^{\otimes \ell}\otimes\boa_r,\ \ \ 1\le r\le N.$$

 We now prove that $\pi$ is injective.  Since $\wt\ker\pi\subset\ell\omega-Q^+$ it follows that there must exist $\bow=(w_1,\cdots, w_N)$, $w_r\in W(\ell,\gr\boa_r)$ with $$\pi(\bow)=0,\ \ \ (x^+\otimes\bc[t])\bow=0.$$ Hence $(x^+\otimes\bc[t])w_r=0$ for all $1\le r\le N$. Using Proposition \ref{ann} and Corollary \ref{hw}, we see that  $w_r\in W(\ell,\gr\boa_r)_\ell$ for all $1\le r\le N$, i.e., we can write $w_r=\boh_rw_{\ell,\gr\boa_r}$, for some $\boh_r\in\bu(\lie h[t])$. This gives,  $$0=\pi(\bow)=\sum_{r=1}^N\pi(w_r)=   \sum_{r=1}^N \boh_r(e_+^{\otimes \ell}\otimes\boa_r) =e_+^{\otimes \ell}\otimes \sum_{r=1}^N\bob_r\boa_r,$$ where $\boh_r(e_+^{\otimes \ell}\otimes 1)= e_+^{\otimes \ell}\otimes\bob_r,$ for some  $\bob_r\in A_\ell^{S_\ell}$,  $1\le r\le N$ by Lemma \ref{impsi}. Theorem \ref{clfree} implies that $\bob_r=0$  and hence $\boh_r(e_+^{\otimes \ell}\otimes 1)= 0$ for  all $1\le r\le N$. Now using Theorem \ref{globalweyl} we conclude that $$w_r=\boh_rw_{\ell, \gr\boa_r}=0,$$ and so $\bow=0$. This proves that $\ker \pi=\{0\}$ and completes the proof of the proposition.
\end{pf}
 Equation \eqref{impsi2} implies that $$\bw(\ell)=\psi_m^k(\bu(\lie g[t])P(m,0)^\ell_\ell). $$
 \begin{cor} There exists an injective  map of $\lie g[t]$--modules $\iota: \bw(\ell)\to P(m,0)^\ell$ such that $\psi_m^k\cdot\iota=\id.$ More explicitly, $\iota(\bw(\ell))$ is a summand of the $\lie g[t]$ submodule of $P(m,0)^\ell$ generated by $P(m,0)^\ell_\ell$ and we have $$\bu(\lie g[t])P(m,0)^\ell_\ell=\iota(\bw(\ell))\oplus \bk,$$ for some $\bk\in\cal I_{\gr}^{\ell\omega}$.
 \end{cor}
 \begin{pf} Since $W(\ell,s)$ is projective in $\cal I_{\gr}^{\ell\omega}$, it follows from  Proposition \ref{direct} that $\bw(\ell)$ is projective. The corollary is now immediate.
 \end{pf}

 \subsection{} Let $\varphi_\ell: P(m,0)^{\ell}\to P(m,0)^{\ell-2}$ be the canonical map of graded modules.  We shall prove, \begin{prop}\label{secondred} For $\ell=m+2k$, we have $$\ker\varphi_\ell=\iota\bw(\ell),$$ and hence we have a short exact sequence,
 $$0\to \bigoplus_{r=1}^NW(\ell ,\gr\boa_r)\to P(m,0)^\ell\to P(m,0)^{\ell-2}\to 0.$$
 \end{prop}

\subsection{}
{\em Proof of Theorem 1.} Assume Proposition \ref{secondred}. Set $$\Gamma_k= \{(m+2i-2)\omega:\;0\leq i\leq k \},\qquad P_k=\Gamma_k P(m,0)=\ker\left(P(m,0)\longrightarrow P(m,0)^{m+2k-2}\right).$$
Proposition \ref{trunc} gives,
$$P_{k}/P_{k+1}=\Gamma_{k}P(m,0)/\Gamma_{k+1} P(m,0)\cong\Gamma_{k}P(m,0)^{m+2k}. $$
It is clear that an equivalent reformulation of Proposition \ref{secondred} is that
$$\Gamma_{k}P(m,0)^{m+2k}\cong  \bigoplus_{r=1}^NW(\ell ,\gr\boa_r).$$
To complete the proof of Theorem 1, we must show that 
\begin{equation}\bigcap_{k\ge 0}P_k=\{0\}. \end{equation} 
It clearly suffices to prove that 
\begin{equation}\label{intersection}P_k[s]=0,\ \ s<k.\end{equation} 
For this, recall from Section \ref{gammatrunc} that 
$$ P_k= \Gamma_k P(m,0)=\sum_{r\ge k}\bu(\lie g[t])P(m,0)_{m+2r}.$$ 
Since $P(m,0)_{m+2r}\subset \bu(\lie g[t])_{r\alpha}p_m,$  and $x^+ p_{m}=0$, we see by using the Poincare--Birkhoff--Witt theorem, that 
$$P(m,0)_{m+2r}\subset\sum_{\bor\in\bn^r}\bu(\lie g[t])\bx_\bor p_m.$$ 
Since $\bx^+_\bor p_m\in P(m,0)[j]$ for $\bor\in\bn^r$ and $j\geq r$, we see that 
$$P(m,0)_{m+2r}\subset \bigoplus_{j\geq r} P(m,0)[j]$$ 
and hence Equation \ref{intersection} is established.

\subsection{} To prove Proposition \ref{secondred}, note that by Corollary \ref{direct} we have,
\begin{equation} \label{pmdirect} P(m,0)^\ell_\ell\ \ =  (\iota \bw(\ell))_\ell\ \ \oplus\ \ \bk_\ell,\end{equation}
as $\lie h$--modules. The proposition follows if we prove that $\bk_\ell=0$.  We begin by  noting that  Lemma \ref{quotwel} implies that  each of the subspaces in \eqref{pmdirect} is  a module for $A_\ell^{S_\ell}$ and hence \eqref{pmdirect} is  an isomorphism of $A_\ell^{S_\ell}$--modules. Moreover we have an isomorphism, $$\iota(\bw(\ell))_\ell\cong\bigoplus_{r=1}^N\left(W(\ell,\gr\boa_r)\right)_\ell \cong\bigoplus_{r=1}^N\tau_r A_\ell^{S_\ell},$$ 
of $A_\ell^{S_\ell}$--modules, where  $\tau_r A_\ell^{S_\ell}$ is the  regular representation of $A_\ell^{S_\ell}$ shifted by grade $\gr\boa_r$. Hence $\iota(\bw(\ell))_\ell$ is a free $A_\ell^{S_\ell}$--module of rank $N$.

Let $\bc_\ell$ be the one dimensional $A_\ell^{S_\ell}$--module obtained by going modulo the augmentation ideal in $A_\ell^{S_\ell}$.
Proposition \ref{secondred} is now clearly a consequence of,
\begin{prop}\label{thirdred}We have $$ \dim_\bc(P(m,0)^\ell_\ell\otimes_{A_\ell^{S_\ell}}\bc_\ell)\le N,\qquad  N=[W_{\loc}(\ell,0): V(m)]=\sum_{s\in\bz_+}[W_{\loc}(\ell, 0): V(m,s)].$$ In particular,
regarded as a right module for $A_\ell^{S_\ell}$, we have that $P(m,0)^\ell_\ell$ is generated by at most $N$ elements.
\end{prop}
Observe that this proposition is clear in the case that $\ell=m$, and so
the rest of the paper is devoted to proving the proposition (and hence Theorem 1) when $\ell > m$. It is helpful to recall here that the  Kostka polynomials involved in the case of $\lie{sl}_2$ are given explicitly by the following,$$K_{m, (m-r\ge r)}(u)= \frac{[m]_u!}{[m-r]_u![r]_u!}-\frac{[m]_u!}{[m-r+1]_u![r-1]_u!},$$ where for an integer $n$ we denote by $$[n]_u=\frac{1-u^n}{1-u},\ \ [n]_u!=[n]_u[n-1]_u\cdots [2]_u,\ \ \ n\in\bz_+,\ \  [n]_u!=0,\ \ n<0. $$ In particular, this means that   $$N=K_{\ell, (\ell-k\ge k)}(1)=\binom{\ell}{k}-\binom{\ell}{k-1}= \sum_{s\in\bz_+}[W_{\loc}(\ell, 0): V(m,s)].$$
\section{Proof of Proposition \ref{thirdred}}\label{proofthirdred}

 We need two key ingredients for the proof: a result which identifies the local Weyl module with a quotient of $\bu(\lie n^-[t])$  established in \cite{CPweyl} and an important identity  established by H. Garland in \cite[Lemma 7.5]{G}.
We use freely the notation of the preceding sections. Thus, $k,m\in\bz_+$ are fixed and we set $\ell=m+2k$. In addition, give  any power series $S(u)$ with coefficients in some algebra $A$, we let $S(u)_r$ be the coefficient of $u^r$ in $S$. Moreover for $a\in A$ we set $$a^{(n)}=\frac{a^n}{n!}.$$ Finally, given any subset $S$ of $A$ we denote by $\langle S\rangle$ the left ideal of $A$ generated by $S$.

\subsection{}
Define formal power series in an indeterminate $u$ with coefficients in $\bu(\lie g[t])$, by
\begin{gather*} X^\pm(u)=\sum_{i\geq 1} (x^\pm\otimes t^i) u^i, \qquad X^\pm_0(u)=\sum_{i\geq 0} (x^\pm\otimes t^i) u^{i+1}, \\ P^+(u) =\exp\left(-\sum_{i\ge 1} \frac{h\otimes t^i}{i} u^i\right). \end{gather*}
Recall  that $\bu(\lie n^\pm[t])$ is the polynomial algebra in the variables $\{x^\pm\otimes t^s: s\in\bz_+\}$ and define ideals  $I(\ell)^\pm$ by \begin{gather*} I(\ell)^\pm = \langle X^\pm_0(u)^{(r)}_s : r\ge 1, s\ge\ell+1\rangle. \end{gather*}
 The adjoint action of $\lie h$ on $\bu(\lie n^\pm[t])$ preserves $I(\ell)^\pm$ and it is simple to see that  \begin{gather*}I(\ell)^\pm_{\pm k\alpha}={\rm{span}}\ \{X^\pm_0(u)^{(r)}_s\bx_\bor^\pm:\bor\in\bz_+^{k-r},\ \ 1\leq r\leq k, s\ge\ell +1\},
\end{gather*} where for any $p\in\bz_+$ and  $\bor=(r_1,\ldots,r_p)\in \bz_+^{p}$,  we recall that 
$$\mathbf x^\pm_{\bor}=(x^\pm\otimes t^{r_1})\cdots (x^\pm\otimes t^{r_p})\in\bu(\lie n^\pm[t]).$$
Regard $\bu(\lie n^-[t])$ as a module for itself via left multiplication and as a $\lie h$--module via the adjoint action. The following is proved in \cite[Lemma 6.4]{CPweyl}.
\begin{prop}\label{idealweyl} The assignment $u\mapsto u(w_{\ell,0}\otimes 1)$ defines a map of $\lie h\oplus \lie n^-[t]$--modules, $\bu(\lie n^-[t])\to W_{\loc}(\ell,0),$ and induces an isomorphism 
$$\bu(\lie n^-[t])/I(\ell)^-\cong W_{\loc}(\ell,0),\qquad \left(\bu(\lie n^-[t])/I(\ell)^-\right)_{-k\alpha}\cong W_{\loc}(\ell,0)_m.$$ 
Moreover, \begin{equation}\label{idealweyl2}\dim_\bc\bu(\lie n^-[t])/I(\ell)^-=2^\ell,\ \ \ \dim_\bc\left(\bu(\lie n^-[t])/I(\ell)^-\right)_{-k\alpha}=\binom{\ell}{k}.\end{equation}\hfill\qedsymbol\end{prop}
  We now prove,
 \begin{cor}\label{mapinj} The map $u\mapsto (x^-\otimes 1)u$ from $\bu(\lie n^-[t])\to\bu(\lie n^-[t])$ induces an injective map of vector spaces $\left(\bu(\lie n^-[t])/I(\ell)^-\right)_{-(k-1)\alpha}\to \left(\bu(\lie n^-[t])/I(\ell)^-\right)_{-k\alpha}.$
 \end{cor}
 \begin{pf} The only statement of the corollary which is not immediate is that the induced map $$\left(\bu(\lie n^-[t])/I(\ell)^-\right)_{-(k-1)\alpha}\to \left(\bu(\lie n^-[t])/I(\ell)^-\right)_{-k\alpha}$$ is injective. To prove this, it suffices to show that the map $(x^-\otimes 1): W_{\loc}(\ell,0)_{m+2}\to W_{\loc}(\ell,0)_{m}$ is injective. But this is clear since $m+2\ge 0$ and $W_{\loc}(\ell, 0)$ is a finite--dimensional $\lie{sl}_2$--module.
 \end{pf}

\subsection{}

The following is a reformulation of a result of H. Garland. In the current form it can be obtained from   \cite[Lemma 1.3]{CPweyl} by applying the anti--automorphism of $\lie g[t]$ which sends $x^+\otimes t^s\mapsto x^-\otimes t^s, h\otimes t^s\mapsto h\otimes t^s$ for all $s\in\bz_+$.
\begin{lem}\label{garland}
Let $s\ge r$.  Then
\begin{eqnarray*}\label{garlandeq} & (x^+\otimes t)^{(s)}(x^-)^{(r)}&=(-1)^r(P^+(u)X^+(u)^{(s-r)})_s + X_{s,r}\\
&&= (-1)^rX(u)^{(s-r)}_s+ H_{s,r}+ X_{s,r},
\end{eqnarray*} 
where \ $X_{s,r}\in \langle \bu(\lie g[t])_{(s-r+1)\alpha}\rangle $\  and \ $H_{s,r}\in (\lie h\otimes t\bc[t])\bu(\lie h\otimes t\bc[t]\oplus\lie  n^+[t])_{(s-r)\alpha}.$  \hfill\qedsymbol
\end{lem}

\subsection{} For the rest of the section we set $$P_{\loc}=P(m,0)^\ell_\ell\otimes_{A_\ell^{S_\ell}}\bc_\ell.$$
\begin{lem} As a vector space $P_{\loc}$ is spanned by the elements $\{\bx_\bor^+( p_{m}\otimes 1):\bor\in\bz_+^k\}.$
\end{lem}
\begin{pf} Recall from Lemma \ref{genpl} that $P(m,0)^\ell_\ell$ is generated as a left $\lie h[t]$--module by the elements $\bx_\bor^+p_m$, $\bor\in\bz_+^k$. Applying Lemma \ref{quotwel} to $P(m,0)^\ell$ we see that $P(m,0)^\ell_\ell$ has the structure of a right $A_\ell^{S_\ell}$--module given by,  
\begin{equation}\label{ploc} v\boa=\boh v,\ \qquad v\in P(m,0)^\ell_\ell,\qquad \boh\in\bu(\lie h[t]),\qquad \boa\in A_\ell^{S_\ell}.\end{equation} Hence $P_{\loc}$ is spanned by the elements $\bx_\bor^+p_m\otimes 1$, $\bor\in\bz_+^k$ as required.

\end{pf}

\subsection{} 
\begin{prop}\label{usegarland}  Let $0\le j\le k$ and set $p=k-j$ and $q\ge m+k+j+1$. Then
\begin{equation}  X^+(u)^{(p)}_q\mathbf x^+_\bor (p_{m}\otimes 1)=0,\ \ {\rm{for\ all}}\ \ \bor\in\bz_+^j.\end{equation}

\end{prop}
\begin{pf} Since $x^+\bx^+_\bor p_{m}=0$ we have $$(x^-)^{m+2j+1}\mathbf x^+_\bor p_{m}=0.$$ Using Lemma \ref{garland} with $s=q$ and $r=q-p$, we get,
 \begin{gather*} 0= (x^+\otimes t)^{(q)}(x^-)^{(m+2j+1)}\mathbf x^+_\bor p_{m}= (-1)^{q-p}X(u)^{(p)}_q\bx^+_\bor p_{m}+ u\bx^+_\bor p_{m}+ u'\bx^+_\bor p_{m},\end{gather*}
 where $u\in(\lie h\otimes t\bc[t])\bu(\lie h\otimes t\bc[t]\oplus\lie n^+[t])_{p\alpha}$ and $u'\in \langle\bu(\lie g[t])_{(p+1)\alpha}\rangle$. Writing 
 $$u=\sum_{s\ge 1}(h\otimes t^s)u_s,\ \ \qquad u_s\in\bu(\lie h\otimes t\bc[t]\oplus\lie n^+[t])_{p\alpha},$$ 
 we get by using Lemma \ref{quotwel} and \eqref{ploc} that,
 $$u\bx_\bor^+(p_{m}\otimes 1)=\sum_{s\ge 1}u_s\bx_\bor^+p_{m}\otimes \symm(h\otimes t^s).1=0.$$
 Further, since $u'\in \langle\bu(\lie g[t])_{(p+1)\alpha}\rangle$ it follows that $u'\bx^+_\bor (p_{m}\otimes 1)=0$ and hence finally, we have $$X^+(u)^{(p)}_q\mathbf x^+_\bor (p_{m}\otimes 1)=0.$$
\end{pf}

Let $U_k$ be the  vector space quotient of $\bu(\lie n^+\otimes t\bc[t])_{k\alpha}$ by the subspace spanned by the elements of the set $\{X^+(u)^{(k-j)}_q\bx^+_\bor:\bor\in\bn^j,\  0\le j\le k,\ \ q\ge m+k+j+1\}$ or, equivalently, the set of elements $ \{X^+(u)^{(p)}_q\bx^+_\bor:\bor\in\bn^{k-p},\ \ 0\le p\le k,\ \ p+q\ge\ell +1\}$. The following is immediate.

\begin{cor} The linear map $ \bu(\lie n^+\otimes t\bc[t])_{k\alpha}\to P_{\loc}$ given by $\bx^+_\bor\mapsto\bx^+_\bor p_{m}$, $\bor\in\bn^k$ is surjective  and induces a surjective map $U_k\to P_{\loc}.$\hfill\qedsymbol

\end{cor}
\subsection{} Using  Corollary \ref{usegarland}, we see that  Proposition \ref{thirdred} is established if we show that \begin{equation}\label{finalred} \dim_\bc U_k\le \binom{\ell}{k}-\binom{\ell}{k-1}.\end{equation} For this, we proceed as follows.
   Set $$J(\ell)^+=\langle X^+(u)^{(p)}_q\ \ : \ \ p+q\ge\ell +1\rangle\subset\bu(\lie n^+[t]).$$ Clearly $J(\ell)^+$ is preserved by the adjoint action of $\lie h$ and we have an isomorphism of vector spaces \begin{equation}\label{jlk} U_k\cong \bu(\lie n^+\otimes t\bc[t])_{k\alpha}/J(\ell)^+_{k\alpha}\cong \left(\bu(\lie n^+\otimes t\bc[t])/J(\ell)^+\right)_{k\alpha}.\end{equation} It is now clear that the following proposition establishes \eqref{finalred}.
  \begin{prop}
   \begin{enumerit}
   \item[(i)] The algebra homomorphism  $\varphi: \bu(\lie n^+[t])\to \bu(\lie n^+\otimes t\bc[t])$ given by sending $$x^+\otimes 1\mapsto 0,\ \qquad  x^+\otimes t^s\mapsto x^+\otimes t^s,\ \ s\ge 1,$$ induces a short exact sequence of algebras $$0\to \langle x^+\otimes 1, \ J(\ell)^+  \rangle / I(\ell)^+\to\bu(\lie n^+[t])/ I(\ell)^+\to\bu(\lie n^+\otimes t\bc[t])/J(\ell)^+ \to 0,$$ and hence also a  short exact sequence of vector spaces,$$0\to \left(\langle x^+\otimes 1,\  J(\ell)^+,\rangle/ I(\ell)^+\right)_{k\alpha}\to\left(\bu(\lie n^+[t])/ I(\ell)^+\right)_{k\alpha}\to\left(\bu(\lie n^+\otimes t\bc[t])/J(\ell)^+\right)_{k\alpha}\to 0.$$
      
       \item[(ii)] We have, $$\dim_\bc\left(\bu(\lie n^+[t])/ I(\ell)^+\right)_{k\alpha}=\binom{\ell}{k},\qquad \dim_\bc\left(\langle x^+\otimes 1,\  J(\ell)^+\rangle/ I(\ell)^+\right)_{k\alpha}\ge\binom{\ell}{k-1}.$$
       \end{enumerit}
  \end{prop}
  \begin{pf}  We first prove that $ I(\ell)^+\subset\langle x\otimes1, J(\ell)^+\rangle$. For this, write $$X_0^+(u)^p_{q+p}=(x^+\otimes1 + X^+(u))^p_{q},$$ and use the binomial expansion to conclude that$$X^+_0(u)^p_{q}\in \langle x\otimes1, J(\ell)^+\rangle,\ \ q+p\ge\ell+1.$$
 Moreover, we also have $$\varphi(X_0^+(u)^p_q)= X^+(u)^p_{q-p}\in J(\ell)^+, \ \ {\rm {if}}\ \ q\ge \ell+1,$$ and hence we have an induced map of algebras $$\bar{\varphi}: \bu(\lie n^+[t])/I(\ell)^+\to\bu(\lie n^+\otimes t\bc[t])/J(\ell)^+\to 0.$$ Suppose that $\varphi(u)\in J(\ell)^+$ for some $u\in\bu(\lie n^+[t])$. Since $\varphi(x\otimes 1)=0$ and the restriction of $\varphi$ to $\bu(\lie n^+\otimes t\bc[t])$ is the identity, we see that $u\in J(\ell)^+$. It follows that the kernel of $\varphi$ is precisely $\langle x\otimes1, J(\ell)^+\rangle$ and the proof is complete.

   To prove (ii) we use Proposition \ref{idealweyl}. Since the map $x^\pm\otimes t^s\to x^\mp\otimes t^s$ induces an isomorphism of algebras,  $$\bu(\lie n^+[t])\cong \bu(\lie n^-[t]),
   \qquad I(\ell)^+\cong I(\ell)^-,$$ we get from \eqref{idealweyl2} that $\dim_\bc\left(\bu(\lie n^+[t])/ I(\ell)^+\right)_{k\alpha}=\binom{\ell}{k}$.

    Consider the algebra homomorphism $\psi: \bu(\lie n^+[t])\to\bu(\lie n^+[t])$ which is given by left multiplication by $(x\otimes 1)$. Clearly $$\psi\left(\bu(\lie n^+[t])_{r\alpha}\right)\subset \bu(\lie n^+[t])_{(r+1)\alpha},\ \  r\in\bz,\qquad \psi(I(\ell)^+)\subset I(\ell)^+,$$ and hence we have a map $$\bar\psi: \bu(\lie n^+[t])/I(\ell)^+\to \bu(\lie n^+[t])/I(\ell)^+,$$ and in fact $$\bar\psi( \bu(\lie n^+[t])/I(\ell)^+)_{(r-1)\alpha}\subset (\langle x\otimes 1, J(\ell)^+\rangle/I(\ell)^+)_{r\alpha}.$$
  By Corollary \ref{mapinj}, the restriction of $\bar\psi$ to $\left(\bu(\lie n^+[t])/I(\ell)^+\right)_{(k-1)\alpha}$ is injective and hence using \eqref{idealweyl2} again, we get
  $$\dim_\bc\left(\langle x^+\otimes 1,\  J(\ell)^+\rangle/ I(\ell)^+\right)_{k\alpha}\ge \dim_\bc(\bu(\lie n^+[t])/I(\ell)^+)_{(k-1)\alpha}=\binom{\ell}{k-1},$$ and part (ii) is proved.

\end{pf}

\end{document}